\documentclass[times,final]{elsarticle}

\usepackage{framed,multirow}
\usepackage{xcolor}
\usepackage[margin=1in]{geometry}

\usepackage{amssymb}
\usepackage{latexsym}

\usepackage{amsmath, amsthm}
\usepackage{multirow}
\usepackage{stmaryrd}
\usepackage{booktabs}
\usepackage{graphicx}
\usepackage{subcaption}
\usepackage{cite}
\usepackage{verbatim}
\usepackage[all,cmtip]{xy}

\renewcommand{\theequation}{\arabic{section}.\arabic{equation}}

\newcommand{\bx}{{\bf x} }

\newcommand{\bJ}{{\bf J} }

\newcommand{\nn}{\nonumber}

\newcommand{\be}{\begin{equation}}
\newcommand{\ee}{\end{equation}}
\newcommand{\ba}{\begin{array}}
\newcommand{\ea}{\end{array}}
\newcommand{\bea}{\begin{eqnarray}}
\newcommand{\eea}{\end{eqnarray}}
\newcommand{\beas}{\begin{eqnarray*}}
\newcommand{\eeas}{\end{eqnarray*}}

\def \qed {\Box}

\newtheorem{theorem}{Theorem}[section]
\newtheorem{lemma}{Lemma}[section]
\newtheorem{rmk}{Remark}[section]

\begin{document}
	
%\verso{Jia Yin}

%\textbf{Highlights}
%\begin{itemize}
%	\item 2
%\end{itemize}  

	\title{A fourth-order compact time-splitting method for the Dirac equation with time-dependent potentials}%
	
	\author[1]{Jia Yin\corref{cor1}}
	\cortext[cor1]{Tel.: +1-9252852235}
    \ead{jiayin@lbl.gov}
	
	\address[1]{Department of Mathematics, National University of Singapore, Singapore 119076, Singapore}

	\begin{abstract}
		%%%
		In this paper, we present an approach to deal with the dynamics of the Dirac equation with time-dependent electromagnetic potentials using the fourth-order compact time-splitting method ($S_\text{4c}$). To this purpose, the time-ordering technique for time-dependent Hamiltonians is introduced, so that the influence of the time-dependence could be limited to certain steps which are easy to treat. Actually, in the case of the Dirac equation, it turns out that only those steps involving potentials need to be amended, and the scheme remains efficient, accurate, as well as easy to implement. Numerical examples in 1D and 2D are given to validate the scheme. 
		%%%%
	\end{abstract}
	

    \maketitle

%\linenumbers


\section{Introduction}\setcounter{equation}{0}
The Dirac equation is a relativistic equation in particle physics which integrates quantum mechanics with special relativity. There has been growing interest in it since it was applied in various areas, such as in graphene and other two-dimensional materials \citep{Graphene, Graphene2, FW, FW2, NGP}, in intense laser field \citep{BCLL, FLB}, in quantum Hall effect \citep{QH, QH2}, and in topological insulators \citep{TI, TI2}.

The Dirac equation with natural units could be represented using the wave function $\Psi:=\Psi(t, \bx)\in\mathbb{C}^4$ in $d$-dimension ($d = 1, 2, 3$) as
\be\label{eq:dirac4}
i\partial_t\Psi =  \left(- i\sum_{j = 1}^{d}\alpha_j\partial_j + \beta\right)\Psi
+ \left(V(t, \mathbf{x})I_4 - \sum_{j = 1}^{d}A_j(t, \mathbf{x})\alpha_j\right)\Psi, \quad t>0, \quad \bx\in\mathbb{R}^d,
\ee
with initial value
\be \label{eq:initial}
\Psi(t=0,\bx)=\Psi_0(\bx),\qquad \bx\in{\mathbb R}^d.
\ee
In the equation, $i$ is the imaginary unit, $t$ represents time, $\bx = (x_1, ..., x_d)^T$ is the spacial coordinate, $\partial_j := \partial_{x_j}$ ($j = 1, ..., d$) are spatial derivatives, and the four-component wave function $\Psi$ could be explicitly written as $\Psi(t, \bx) = (\psi_1(t, \bx), \psi_2(t, \bx), \psi_3(t, \bx), \psi_4(t, \bx))^T$. $V(t, \bx)$ and $\mathbf(t, \bx) := (A_1(t, \bx), ..., A_d(t, \bx))^T$ are real functions, which serve as the electric and the magnetic potentials, respectively. Moreover, $I_n$ is the $n\times n$ identity matrix, while $\alpha_j$ ($j = 1, ..., d$) and $\beta$ are $4\times 4$ Dirac matrices defined as
\be 
\begin{aligned}
	\label{eq:alpha}
	\alpha_1=\left(\begin{array}{cc}
		\mathbf{0} & \sigma_1  \\
		\sigma_1 & \mathbf{0}  \\
	\end{array}
	\right),\quad
	\alpha_2=\left(\begin{array}{cc}
		\mathbf{0} & \sigma_2 \\
		\sigma_2 & \mathbf{0} \\
	\end{array}
	\right), \quad
	\alpha_3=\left(\begin{array}{cc}
		\mathbf{0} & \sigma_3 \\
		\sigma_3 & \mathbf{0} \\
	\end{array}
	\right),\quad
	\beta=\left(\begin{array}{cc}
		I_{2}& \mathbf{0} \\
		\mathbf{0} & -I_{2} \\
	\end{array}
	\right),
\end{aligned}
\ee
with the Pauli matrices
\be\label{eq:Paulim}
\sigma_{1}=\left(
\begin{array}{cc}
	0 & 1  \\
	1 & 0  \\
\end{array}
\right), \qquad
\sigma_{2}=\left(
\begin{array}{cc}
	0 & -i \\
	i & 0 \\
\end{array}
\right),\qquad
\sigma_{3}=\left(
\begin{array}{cc}
	1 & 0 \\
	0 & -1 \\
\end{array}
\right).
\ee

The dynamics of the Dirac equation \eqref{eq:dirac4} has been widely studied both analytically and numerically. The dispersion relation suggests that the wavelength is at $O(1)$ in space and time. For the existence and multiplicity of bound states and/or standing wave solutions, we refer to \citep{Das,DK,ES,GGT,Gross,Ring} and references therein. On the other hand, many efficient and accurate numerical methods have been proposed and analyzed \citep{AL,BL}, such as the finite difference time domain (FDTD) methods \citep{ALSFB,HP,NSG}, splitting methods \citep{BCJT2,BSG,FLB0,HJM,MK}, exponential wave integrator Fourier pseudospectral (EWI-FP) method \citep{BCJT2}, the Gaussian beam method \citep{WHJY}, etc. For atomic processes in relativistic heavy-ion collisions, a treatment in momentum space was introduced in \citep{MBS}. Additionally, there have been many studies on different regimes of the Dirac equation, such as the nonrelativistic regime \citep{BCJT, BCJY, BCY, CaiW2}, and the semiclassical regime \citep{MY}.

In order to increase the convergence rate in time while maintain a relatively small computational cost, a fourth-order compact time-splitting method ($S_\text{4c}$) was introduced for the Dirac equation \citep{BY}. Compared to other fourth-order splitting methods, such as the Forest-Ruth scheme ($S_4$) \citep{FR} (for the Dirac equation, $S_4$ has been applied in \citep{BSG}), and the partitioned Runge-Kutta scheme ($S_\text{4RK}$) \citep{Geng}, $S_\text{4c}$ is more efficient, and avoids negative time steps in sub-problems.  However, the method in \citep{BY} is only valid for time-independent potentials, i.e., $V(t, \bx) \equiv V(\bx)$, $\mathbf{A}(t, \bx) \equiv \mathbf{A}(\bx)$ in \eqref{eq:dirac4}. When the potentials are time-dependent, it is not straightforward to extend the method, resulting in the limitation in application. In this paper, we apply the time ordering technique, which was introduced in \citep{CC2}, so that the extension to time-dependent potentials could be realized. Numerical tests are also carried out to validate the extension of $S_\text{4c}$, and compare its performance with other methods.

For simplicity, the majority of this paper only deals with the Dirac equation in one dimension (1D) and two dimensions (2D). As given in \citep{BCJT2}, in 1D and 2D, \eqref{eq:dirac4} could be reduced to 	
\be
\label{eq:dirac2}
i\partial_t\Phi = \left(-i\sum_{j = 1}^{d}\sigma_j\partial_j + \sigma_3\right)\Phi + \left(V(t, \mathbf{x})I_2 - \sum_{j = 1}^{d}A_j(t, \mathbf{x})\sigma_j\right)\Phi, \quad \bx \in\mathbb{R}^d, \quad d = 1, 2,
\ee
with the initial condition
\be\label{eq:initial11}
\Phi(t = 0, \mathbf{x}) = \Phi_0(\mathbf{x}), \quad \mathbf{x}\in\mathbb{R}^d, \quad d = 1, 2,
\ee
where the two-component wave function $\Phi = (\psi_1, \psi_4)^T$ (or $\Phi=(\psi_2,\psi_3)^T$). Extension of the results to the four-component equation \eqref{eq:dirac4} is straightforward. 

The rest of the paper is organized as follows. In section 2, a review of the time-ordering technique for time-dependent Hamiltonians is given. The application of the technique to $S_\text{4c}$ for the Dirac equation with time-dependent electromagnetic potentials is discussed in section 3. Section 4 shows numerical results in 1D and 2D to numerically validate the scheme, and conclusions are drawn in section 5.\\

\section{The time-ordering technique for time-dependent Hamiltonians}
The main idea to deal with the time-dependent potentials in the Dirac equation when applying splitting methods is to use the time-ordering technique, which was first introduced by Suzuki in \citep{Suzuki2}. The idea has been successfully applied to the Schr\"{o}dinger equation with time-dependent potentials \citep{CC2}. For the Dirac equation, time-ordering for the splitting method was mentioned in \citep{FLB0}. In that paper, time-ordering was omitted in the end because the error introduced is second-order in time step, which is the same as the order of the splitting method there.

In this section, we give a detailed explanation of the time-ordering technique, where the key point is given in Lemma \ref{lemma:to}.

For illustration, we first consider a model equation ($d = 1, 2, 3$)
\begin{equation}
\label{eq:model}
\partial_tu(t, \bx) = (T + W(t))u(t, \bx), \quad t>t_0, \quad \bx\in\mathbb{R}^d,
\end{equation}
with the initial data
\begin{equation}
\label{eq:it55}
u(t_0, \bx) = u_0(\bx), \quad \bx\in\mathbb{R}^d,
\end{equation}
where $t_0$ is the initial time, $T$ is a time-independent operator, and $W(t)$ is a time-dependent one. We remark here that the wave function $u(t, \bx)$ could either be a scalar or a vector function. Here we focus on the temporal coordinate, so the spatial coordinates are not taken into account in the expression of the operators, and we can further take  $u(t) := u(t, \bx)$ for simplicity. Denote $H(t) := T + W(t)$, suppose the exact solution $u(t)$ propagates with the operator $U(t, t_0)$, i.e.,
\begin{equation}\label{eq:iter}
u(t) = U(t, t_0)u(t_0),
\end{equation}
then by plugging \eqref{eq:iter} into \eqref{eq:model}, we can get the differential equation
\begin{equation}\label{eq:model_U}
\partial_tU(t, t_0) = H(t)U(t, t_0), \quad t>t_0,
\end{equation}
with $U(t_0, t_0) = Id$, the identity operator, which can easily be checked. Take any $\tau>0$, then by Taylor expansion,
\be\label{eq:Taylor}
\begin{aligned}
	U(t_0 + \tau, t_0) &= U(t_0, t_0) + \tau\partial_tU(t_0, t_0) + O(\tau^2)\\
	&= \left(Id + \tau H(t_0)\right) + O(\tau^2) = e^{\tau H(t_0)} + O(\tau^2).
\end{aligned}
\ee
Noticing the fact that
\begin{equation}
U(t+\tau, t) = \Pi_{k=1}^n U\left(t+\frac{k}{n}\tau, t + \frac{k-1}{n}\tau\right) = \Pi_{k=1}^ne^{\frac{\tau}{n}H\left(t+\frac{k-1}{n}\tau\right)} + O\left(\frac{\tau^2}{n^2}\right),
\end{equation}
where the relation \eqref{eq:Taylor} is used to get the second equality, holds for any positive integer $n$, we have
\begin{equation}\label{eq:repre1}
U(t+\tau, t) = \lim_{n\rightarrow\infty}e^{\frac{\tau}{n}H\left(t+\frac{k-1}{n}\tau\right)}...e^{\frac{\tau}{n}H\left(t+\frac{\tau}{n}\right)}e^{\frac{\tau}{n}H\left(t\right)}, \quad t>0.
\end{equation}

On the other hand, from \eqref{eq:model_U} with the initial condition $U(t_0, t_0) = Id$, we have
\begin{eqnarray}
U(t, t_0) &=& Id + \int_{t_0}^t H(s)U(s, t_0)ds\nonumber\\
&=& Id + \int_{t_0}^tH(s_1)ds_1 + \int_{t_0}^tH(s_1)\int_{t_0}^sH(s_2)U(s_2, t_0)ds_2ds_1\nonumber\\
&=& Id + \sum_{n=1}^\infty\int_{t_0}^t\int_{t_0}^{s_1}...\int_{t_0}^{s_{n-1}}ds_{n}...ds_1H(s_1)...H(s_n)\label{eq:torder}\\
&=:& \mathcal{T}(e^{\int_{t_0}^tH(s)ds}),
\end{eqnarray}
where $\mathcal{T}(\cdot)$ is defined as \textbf{the time-ordering operator}, with the expression given in \eqref{eq:torder}. This gives us
\be\label{eq:repre2}
U(t+\tau, t) = \mathcal{T}\left(e^{\int_{t}^{t+\tau}H(s)ds}\right), \quad t>0.
\ee

From the above discussion, combining \eqref{eq:repre1} and \eqref{eq:repre2}, we get
\be
\mathcal{T}\left(e^{\int_{t}^{t+\tau}H(s)ds}\right)= \lim_{n\rightarrow\infty}e^{\frac{\tau}{n}H\left(t+\frac{k-1}{n}\tau\right)}...e^{\frac{\tau}{n}H\left(t+\frac{\tau}{n}\right)}e^{\frac{\tau}{n}H(t)}, \quad t>0.
\ee

Define a forward time derivative operator \citep{CC2} $\mathcal{D}:= \frac{\stackrel{\leftarrow}{\partial}}{\partial t}$, which is applied to the function on the left-hand side, by
\begin{equation}
	f(t)\mathcal{D} = \lim\limits_{\tau\rightarrow 0}\frac{f(t+\tau) - f(t)}{\tau},
\end{equation}
for any time-dependent function $f(t)$. It is straight forward to observe that
\begin{equation}\label{eq:td_operator}
F(t)e^{\tau\mathcal{D}}G(t) = F(t+\tau)G(t), \quad t>0,
\end{equation}
where $F(\cdot)$ and $G(\cdot)$ are any two time-dependent operators. Then we have the following important lemma.\\

\begin{lemma} 
	\label{lemma:to}
	The following equality holds true for any time-dependent operator $H(t)$.
	\begin{equation}\label{eq:torder_e}
	\mathcal{T}\left(e^{\int_t^{t+\tau}H(s)ds}\right) = \exp[\tau(H(t)+\mathcal{D})], \quad t>0.
	\end{equation}
\end{lemma}
\proof{
	We start from the right-hand-side of \eqref{eq:torder_e}.
	\begin{eqnarray}
	\exp[\tau(H(t)+\mathcal{D})] &=& \lim_{n\rightarrow\infty}\left(e^{\frac{\tau}{n}\mathcal{D}}e^{\frac{\tau}{n}H(t)}\right)^n\nonumber\\
	&=& \lim_{n\rightarrow\infty} e^{\frac{\tau}{n}\mathcal{D}}e^{\frac{\tau}{n}H(t)}...e^{\frac{\tau}{n}\mathcal{D}}e^{\frac{\tau}{n}H(t)}e^{\frac{\tau}{n}\mathcal{D}}e^{\frac{\tau}{n}H(t)}\nonumber\\
	&=& \lim_{n\rightarrow\infty}e^{\frac{\tau}{n}H\left(t+\frac{k-1}{n}\tau\right)}...e^{\frac{\tau}{n}H\left(t+\frac{\tau}{n}\right)}e^{\frac{\tau}{n}H(t)}\nonumber\\
	&=& \mathcal{T}\left(e^{\int_{t}^{t+\tau}H(s)ds}\right),
	\end{eqnarray}
	which gives us the expected result. The first equality in the proof comes from the fact that $e^{x(A+B)} = \lim_{n\rightarrow\infty}\left(e^{\frac{x}{n}A}e^{\frac{x}{n}B}\right)^n$ \citep{CC2}. $\qed$
}

Recall $H(t) = T + W(t)$ in the model equation \eqref{eq:model}. Define $\widetilde{T} = T + \mathcal{D}$, then from \eqref{eq:repre2} and the above lemma, $u(t+\tau)$ can be expressed as
\begin{eqnarray}
u(t+\tau) &=& U(t+\tau, t)u(t) = \mathcal{T}\left(e^{\int_t^{t+\tau}\left(T + W(s)\right)ds}\right)u(t) = \exp[\tau(T + W(t) + \mathcal{D})]u(t)\nonumber\\
&=& \exp[\tau(\widetilde{T} + W(t))]u(t), \quad t>0, \label{eq:td_exact}
\end{eqnarray}
which serves as the exponential expression of the exact solution to \eqref{eq:model}.\\

\section{$S_\text{4c}$ for the Dirac equation with time-dependent potentials}
In this section, $S_\text{4c}$ is applied to the Dirac equation in 1D and 2D. The application is then generalized to the Dirac equation in 3D. Mass conservation and convergence of the method are presented in the last subsection.

\subsection{$S_\text{4c}$ in 1D and 2D}
Based on the time-ordering technique introduced in the previous section, we can now apply $S_\text{4c}$ to the Dirac equation \eqref{eq:dirac2} with time-dependent electromagnetic potentials $V(t, \bx)$ and $\mathbf{A}(t, \bx)$.

Define
\be\label{eq:split_TW}
T = -\sum_{j=1}^d\sigma_j\partial_j - i\sigma_3, \quad
W(t) = -i\left(V(t, \bx)I_2 - \sum_{j=1}^dA_j(t, \bx)\sigma_j\right), \quad d = 1, 2,
\ee
then the Dirac equation \eqref{eq:dirac2} can be expressed in the form \eqref{eq:model} with $u(t, \bx) := \Phi(t, \bx)$.

Applying $S_\text{4c}$ \citep{BY, Chin, CC, CC2} to the exact solution \eqref{eq:td_exact} with time step size $\tau$, we get
\begin{eqnarray}\label{eq:S4c}
\Phi(t+\tau)\approx S_\text{4c}(\tau)\Phi(t) &:=& e^{\frac{\tau}{6}W(t)}e^{\frac{\tau}{2}\widetilde{T}}e^{\frac{2\tau}{3}\widehat{W}(t)}e^{\frac{\tau}{2}\widetilde{T}}e^{\frac{\tau}{6}W(t)}\Phi(t)\\
&=& e^{\frac{\tau}{6}W(t+\tau)}e^{\frac{\tau}{2}T}e^{\frac{2\tau}{3}\widehat{W}\left(t+\frac{\tau}{2}\right)}e^{\frac{\tau}{2}T}e^{\frac{\tau}{6}W(t)}\Phi(t),\nonumber
\end{eqnarray}
where the relation in \eqref{eq:td_operator} is used. In the expression, we have
\begin{eqnarray}
\widehat{W}(t) &:=& W(t) + \frac{\tau^2}{48}[W(t), [\widetilde{T}, W(t)]]\nonumber\\
&=& W(t) + \frac{\tau^2}{48}[W(t), [T, W(t)]] + \frac{\tau^2}{48}[W(t), [\mathcal{D}, W(t)]],
\end{eqnarray}
Through simple computation, we can obtain
\begin{eqnarray}
[W(t), [\mathcal{D}, W(t)]] &=& [W(t), [\mathcal{D}W(t) - W(t)\mathcal{D}]] =  [W(t), [\mathcal{D}W(t) - (W'(t) + \mathcal{D}W(t))]]\nonumber\\
&=& [W(t), -W'(t)] = 0.
\end{eqnarray}
As a result, 
\begin{equation}\label{eq:W_hat}
\widehat{W}(t) = W(t) + \frac{\tau^2}{48}[W(t), [T, W(t)]],
\end{equation}
and the double commutator $[W(t), [T, W(t)]]$ could be represented as shown in the following lemmas.\\

\begin{lemma}\label{lemma:db_com_1d}
	The explicit form of the double commutator $[W(t), [T, W(t)]]$ for the Dirac equation \eqref{eq:dirac2} in 1D ($d=1$) with the splitting \eqref{eq:split_TW} is
	\be\label{eq:commu_1d}
	[W(t), [T, W(t)]]  = -4iA_1^2(t, x)\sigma_3.
	\ee
\end{lemma}
The details of the derivation could be found in Appendix A. 

Similar to the 1D case, we can obtain the double commutator in 2D ($d=2$):
\begin{lemma}\label{lemma:db_com_2d}
	The explicit form of the double commutator $[W(t), [T, W(t)]]$ for the Dirac equation \eqref{eq:dirac2} in 2D ($d=2$) with the splitting \eqref{eq:split_TW} is
    \be\label{eq:commu_2d}
    [W(t), [T, W(t)]] = F_3(t, \mathbf{x}) + F_1(t, \mathbf{x})\partial_1 + F_2(t, \mathbf{x})\partial_2, \quad \text{when }d=2,
    \ee
    where
    \beas
    F_1(t, \mathbf{x}) &=& 4\Bigl( - A_2^2(t, \bx)\sigma_1+A_1(t, \bx)A_2(t, \bx)\sigma_2\Bigr),
    \quad F_2(t, \mathbf{x}) = 4\Bigl(A_1(t, \bx)A_2(t, \bx)\sigma_1 - A_1^2(t, \bx)\sigma_2\Bigr), \\
    F_3(t, \mathbf{x}) &=& 4
    \Big(A_1(t, \bx)\partial_2A_2(t, \bx)-A_2(t, \bx)\partial_1A_2(t, \bx)\Big)
    \sigma_1+4\Big(A_2(t, \bx)\partial_1A_1(t, \bx)-A_1(t, \bx)\partial_2A_1(t, \bx)\Big)
    \sigma_2\\
    &&+4i\left(A_2(t, \bx)\partial_1V(t, \bx)-A_1(t, \bx)\partial_2 V(t, \bx)-
    \big(A_1^2(t, \bx) + A_2^2(t, \bx)\big)\right)\sigma_3.
    \eeas
\end{lemma}
The details of the derivation could be found in Appendix B.\\

From Lemmas \ref{lemma:db_com_1d} and \ref{lemma:db_com_2d}, noticing \eqref{eq:S4c}, the semi-discretized fourth-order compact time-splitting method ($S_\text{4c}$) for the Dirac equation \eqref{eq:dirac2} in 1D and 2D with time-dependent electromagnetic potentials could be defined as:
\be\label{eq:S4c_td}
\Phi^{n+1}(\bx) = e^{\frac{1}{6}\tau W(t_{n+1})}e^{\frac{1}{2}\tau T}e^{\frac{2}{3}\tau\widehat{W}(t_n+\tau/2)}
e^{\frac{1}{2}\tau T}e^{\frac{1}{6}\tau W(t_n)}\Phi^n(\bx), \quad 0\leq n\leq \frac{T}{\tau}-1,
\ee
with the given initial value
\be
\Phi^0(\bx) := \Phi_0(\bx), \quad \bx\in\mathbb{R}^d, \quad d = 1, 2.
\ee
The solution is computed until $T_\text{max}>0$. In the scheme, $\widehat{W}(t)$ is defined as \eqref{eq:W_hat} with $[W(t), [T, W(t)]]$ given in \eqref{eq:commu_1d} and \eqref{eq:commu_2d} respectively for 1D and 2D cases.
$\Phi^n(\bx)$ is the semi-discretized approximation of $\Phi(t, \bx)$ at $t = t_n := n\tau$.

\begin{rmk}
	The application of $S_\text{4c}$ in 1D and 2D to \eqref{eq:dirac2} can be easily extended to the four-component Dirac equation \eqref{eq:dirac4}. Similar to the two-component case, we get
	\begin{equation}
	\Psi(t+\tau)\approx S_\text{4c}(\tau)\Psi(t) = e^{\frac{\tau}{6}W(t+\tau)}e^{\frac{\tau}{2}T}e^{\frac{2\tau}{3}\widehat{W}\left(t+\frac{\tau}{2}\right)}e^{\frac{\tau}{2}T}e^{\frac{\tau}{6}W(t)}\Psi(t),
	\end{equation}
	where
	\begin{equation}
	\widehat{W}(t) = W(t) + \frac{\tau^2}{48}[W(t), [T, W(t)]].
	\end{equation}
	
	For the four-component Dirac equation \eqref{eq:dirac4} in 1D, under the splitting 
	\begin{equation}
	T = -\alpha_1\partial_1 - i\beta, \qquad
	W = -i\Bigl(V(t, x)I_4 - A_1(t, x)\alpha_1\Bigr),
	\end{equation}
	the double commutator $[W(t), [T, W(t)]]$  could be easily derived as:
	\be\label{eq:commut_1d_4}
	[W(t), [T, W(t)]]  = -4iA_1^2(t, x)\beta.
	\ee
	
	For the four-component Dirac equation \eqref{eq:dirac4} in 2D, under the splitting 
	\begin{equation}
	\label{TW2d2}
	T = -\alpha_1\partial_1 -
	\alpha_2\partial_2 - i\beta,
	\quad W = -i\Bigl(V(t, \bx)I_2 - A_1(t, \bx)\alpha_1 - A_2(t, \bx)\alpha_2\Bigr),
	\end{equation}
	the double commutator $[W(t), [T, W(t)]]$  could be easily derived as:
	\be\label{eq:commut_2d_4}
	[W(t), [T, W(t)]] = F_3(t, \mathbf{x})+ F_1(t, \mathbf{x})\partial_1 + F_2(t, \mathbf{x})\partial_2,
	\ee
	where
	\beas
	%\label{commutator_2D2_F}
	F_1(t, \mathbf{x}) &=& 4\Bigl( - A_2^2(t, \bx)\alpha_1+A_1(t, \bx)A_2(t, \bx)\alpha_2\Bigr),
	\quad F_2(t, \mathbf{x}) = 4\Bigl(A_1(t, \bx)A_2(t, \bx)\alpha_1 - A_1^2(t, \bx)\alpha_2\Bigr), \\
	F_3(t, \mathbf{x}) &=& 4\Big(A_1(t, \bx)\partial_2A_2(t, \bx)-
	A_2(t, \bx)\partial_1A_2(t, \bx)\Big)\alpha_1+
	4\Big(A_2(t, \bx)\partial_1A_1(t, \bx)-
	A_1(t, \bx)\partial_2A_1(t, \bx)\Big)\alpha_2\\
	&&+4i\Big(A_2(t, \bx)\partial_1V(t, \bx)-
	A_1(t, \bx)\partial_2V(t, \bx)\Big)\gamma\alpha_3
	- 4i
	\big(A_1^2(t, \bx) + A_2^2(t, \bx)\big)\beta,
	\eeas
	with
	\be\label{gammam}
	\gamma = \begin{pmatrix}  \mathbf{0} &I_2\\ I_2 &\mathbf{0} \end{pmatrix}. \quad
	\ee
\end{rmk}

From the above remark, noticing \eqref{eq:S4c}, the semi-discretized fourth-order compact time-splitting method ($S_\text{4c}$) for the Dirac equation \eqref{eq:dirac4} in 1D and 2D with time-dependent electromagnetic potentials could be defined as:
\be\label{eq:S4c_4}
\Psi^{n+1}(\bx) = e^{\frac{1}{6}\tau W(t_{n+1})}e^{\frac{1}{2}\tau T}e^{\frac{2}{3}\tau\widehat{W}(t_n+\tau/2)}
e^{\frac{1}{2}\tau T}e^{\frac{1}{6}\tau W(t_n)}\Psi^n(\bx), \quad 0\leq n\leq \frac{T}{\tau}-1,
\ee
with the given initial value
\be\label{eq:ini}
\Psi^0(\bx) := \Psi_0(\bx), \quad \bx\in\mathbb{R}^d, \quad d = 1, 2.
\ee
The solution is computed until $T_\text{max}>0$. In the scheme, $\widehat{W}(t)$ is defined as \eqref{eq:W_hat} with $[W(t), [T, W(t)]]$ given in \eqref{eq:commut_1d_4} and \eqref{eq:commut_2d_4} respectively for 1D and 2D cases.
$\Psi^n(\bx)$ is the semi-discretized approximation of $\Psi(t, \bx)$ at $t = t_n := n\tau$.

\subsection{$S_\text{4c}$ in 3D}
In the 3D case, we consider the four-component Dirac equation \eqref{eq:dirac4}. The following lemma shows the application of $S_\text{4c}$ in 3D:
\begin{lemma}
	\label{DTW3d}
	For the Dirac equation \eqref{eq:dirac4} in 3D, i.e. $d = 3$, define
	\begin{equation}
	T = -\sum_{j = 1}^3\alpha_j\partial_j -i\beta,\quad
	W(t) = -i\Bigl(V(t, \bx)I_4 - \sum_{j = 1}^3A_j(t, \bx)\alpha_j\Bigr),
	\end{equation}
	we have
	\begin{equation}
	\label{commun_3d}
	[W(t), [T, W(t)]] = F_4(t, \mathbf{x})+F_1(t, \mathbf{x})\partial_1 + F_2(t, \mathbf{x})\partial_2 + F_3(t, \mathbf{x})\partial_3,
	\end{equation}
	where
	\beas
	F_1(t, \mathbf{x}) &=&4\Big(-\big(A_2^2(t) + A_3^2(t)\big)\alpha_1 +
	A_1(t)A_2(t)\alpha_2 + A_1(t)A_3(t)\alpha_3\Big), \\
	F_2(t, \mathbf{x}) &=&4\Big(
	A_2(t)A_1(t)\alpha_1-\big(A_1^2(t) + A_3^2(t)\big)\alpha_2 + A_2(t)A_3(t)\alpha_3\Big), \\
	F_3(t, \mathbf{x}) &=&4\Big(
	A_3(t)A_1(t)\alpha_1 + A_3(t)A_2(t)\alpha_2-\big(A_1^2(t) + A_2^2(t)\big)\alpha_3\Big), \\
	F_4(t, \mathbf{x}) &=&4
	\Big(A_1(t)\big(\partial_2A_2(t)+\partial_3A_3(t)\big)-
	A_2(t)\partial_1A_2(t)-A_3(t)\partial_1A_3(t)
	\Big)\alpha_1\\
	&&+4
	\Big(A_2(t)\big(\partial_1A_1(t)+\partial_3A_3(t)\big)-
	A_1(t)\partial_2A_1(t)-A_3(t)\partial_2A_3(t)
	\Big)\alpha_2\\
	&&+4
	\Big(A_3(t)\big(\partial_1A_1(t)+\partial_2A_2(t)\big)-
	A_1(t)\partial_3A_1(t)-A_2(t)\partial_3A_2(t)
	\Big)\alpha_3\\
	&&+4i
	\Big(A_1(t)\big(\partial_2A_3(t)-\partial_3A_2(t)\big)+
	A_2(t)\big(\partial_3A_1(t)-\partial_1A_3(t)\big)\\
	&&+ A_3(t)\big(\partial_1A_2(t)-\partial_2A_1(t)\big)
	\Big)\gamma + 4i\Big(A_3(t)\partial_2V(t) - A_2(t)\partial_3V(t)\Big)\gamma\alpha_1\\
	&& + 4i\Big(A_1(t)\partial_3V(t) - A_3(t)\partial_1V(t)\Big)\gamma\alpha_2\\
	&&+ 4i\Big(A_2(t)\partial_1V(t) - A_1(t)\partial_2V(t)\Big)\gamma\alpha_3
	- 4i\Big(A_1^2(t) + A_2^2(t) + A_3^2(t)\Big)\beta.
	\eeas
	Here we use $V(t) := V(t, \bx)$ and $A_j(t) := A_j(t, \bx)$, $j = 1, 2, 3$ for simplicity.
\end{lemma}
The details of the proof could be found in Appendix C.\\

From Lemma \ref{DTW3d}, the semi-discretized fourth-order compact time-splitting method ($S_\text{4c}$) for the Dirac equation \eqref{eq:dirac4} in 3D with time-dependent electromagnetic potentials could be defined in the same way as \eqref{eq:S4c_4} with the initial value \eqref{eq:ini}. Under this circumstance, $\widehat{W}(t)$ is defined as \eqref{eq:W_hat} with $[W(t), [T, W(t)]]$ given in \eqref{commun_3d}.\\
 
According to the explicit forms of the double commutators, we could see that their existence is closely related to the magnetic potentials. In other words, as long as $A_j(t, \bx) \equiv 0$, $j = 1, ..., d$, $\widehat{W}(t) \equiv W(t)$, and the step involving $\widehat{W}(t)$ will have no difference with the steps of $W(t)$. If $A_j(t, \bx) \not\equiv 0$ for some $j = 1, ..., d$, then in 1D, it is still straightforward to compute, but in 2D or 3D, the step involving $\widehat{W}$ will be much more difficult to deal with. Similar to the discussions in \citep{BY}, we may use the method of characteristics and the nonuniform fast Fourier transform (NUFFT) to evaluate the operator involving $\widehat{W}$.

We remark here that similar to other splitting methods, this method could be efficiently applied to different regimes of the Dirac equation. Details are omitted here for brevity.

\subsection{Mass conservation and convergence}
$S_\text{4c}$ with time-dependent potentials conserves mass, as shown in the following lemma.
\begin{lemma}
	For any $\tau>0$, the $S_\text{4c}$ method \eqref{eq:S4c_td} for \eqref{eq:dirac2} conserves the mass, i.e., for $d = 1, 2$
	\begin{equation} \label{eq:mass2}
	\Big\|\Phi^{n+1}\Big\|^2_{L^2} := \int_{\mathbb{R}^d}\left|\Phi^{n+1}\right|^2 d\mathbf{x} = \int_{\mathbb{R}^d}\left|\Phi^{0}\right|^2 d\mathbf{x} = \int_{\mathbb{R}^d}\left|\Phi_0\right|^2 d\mathbf{x} = \Big\|\Phi_0\Big\|^2_{L^2},\quad n\ge0.
	\end{equation}
	Mass conservation also holds for \eqref{eq:S4c_4} to solve \eqref{eq:dirac4}, i.e., for $d = 1, 2, 3$
	\begin{equation} \label{eq:mass4}
	\Big\|\Psi^{n+1}\Big\|^2_{L^2} := \int_{\mathbb{R}^d}\left|\Psi^{n+1}\right|^2 d\mathbf{x} = \int_{\mathbb{R}^d}\left|\Psi^{0}\right|^2 d\mathbf{x} = \int_{\mathbb{R}^d}\left|\Psi_0\right|^2 d\mathbf{x} = \Big\|\Psi_0\Big\|^2_{L^2},\quad n\ge0.
	\end{equation}
\end{lemma}
Proof of the lemma is similar to the proof in \citep{BY}. The details are omitted here for brevity.\\

Moreover, for any $T_\text{max}>0$, define the error function
\begin{equation}
\mathbf{e}^n(\bx) = \Phi(t_n, \bx) - \Phi^n(\bx), \quad 0\leq n\leq \frac{T_\text{max}}{\tau}, \quad \bx\in\mathbb{R}^d, \quad d = 1, 2
\end{equation}
for \eqref{eq:dirac2}, and 
\begin{equation}
\mathbf{e}^n(\bx) = \Psi(t_n, \bx) - \Psi^n(\bx), \quad 0\leq n\leq \frac{T_\text{max}}{\tau}, \quad \bx\in\mathbb{R}^d, \quad d = 1, 2, 3
\end{equation}
for \eqref{eq:dirac4}, then the error bound for $S_\text{4c}$ is given in Theorem \ref{thm:err}.

\begin{theorem}\label{thm:err}
	Let $\Phi^n(\bx)$ be the numerical approximation obtained from $S_\text{4c}$ \eqref{eq:S4c_td} for \eqref{eq:dirac2} (or \eqref{eq:S4c_4} for \eqref{eq:dirac4}),
	then under certain regularity conditions,
	we have the following error estimate
	\be
	\|{\bf e}^n(x)\|_{L^2}\lesssim \tau^4, \quad 0\le n\le\frac{T}{\tau}.
	\ee
\end{theorem}
The idea of the proof is similar to the proof in \citep{BCJY}, so for brevity, the details are omitted here.

\section{Numerical results}
This section consists of numerical examples in 1D and 2D to verify the accuracy of $S_\text{4c}$ \eqref{eq:S4c_td} for the Dirac equation with time-dependent electromagnetic potentials. 

\subsection{Klein paradox}
We first consider a special phenomenon for the Dirac equation, which is called the `Klein paradox' \citep{BSG,FLB0,KSG}, to validate our algorithm. `Klein paradox' describes the different reflection and transmission behavior of the Dirac equation from those of the non-relativistic Schr\"{o}dinger equation of the plane wave solution under a step potential \citep{Klein}. 

Suppose we have a step potential with heigt $V_0$. In the Schr\"{o}dinger case, when the wave energy $E<V_0$, the transmission coefficient is very small, which means most of the wave function is reflected. By contrast, in the Dirac case, when $E<V_0-mc^2$, there could be a non-negligible transmission coefficient. It is believed that the transmitted part comes from the negative energy solution for anti-fermions, while the reflected part is related to the solution for fermions \citep{DC,Greiner,GMR,KSG}. This numerical test is chosen here because there is an analytical transmission coefficient, so that we could compare it with our numerical results.

\begin{figure}[htbp]
	\centering
	\includegraphics[width=0.7\textwidth]{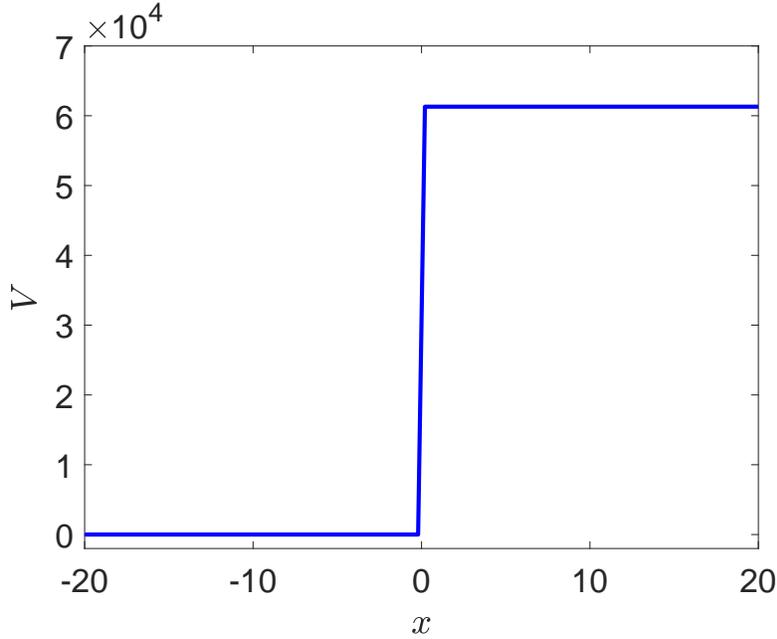}
	\caption{$V(x)$ with $L=10^{-4}$ and $V_0 = 6.13\times 10^4$.}
	\label{fig:KP_V}
\end{figure}

In this example, we consider the 1D Dirac equation
\begin{equation}
	i\partial_t\Phi(t, x) = \left(-ic\sigma_1\partial_x + mc^2\sigma_3\right)\Phi(t, x) + e\left(V(t, x)I_2 - A_1(t, x)\sigma_1\right)\Phi(t, x), \quad t>0, \; x\in\mathbb{R}, 
\end{equation}
where $\Phi(t,x)$ is a two-component wave function, $c$ is the light velocity, $m$ is the fermion mass, and $e$ refers to the electric charge. Specifically, here we take the atomic units, where $c=1/\alpha$ with $\alpha$ being the fine structure constant $\alpha\approx 1/137.0359895$, $m=1$ and $e = 1$. The magnetic potential is taken to be $0$, and the electric potential is given by
\begin{equation}\label{eq:KP_potential}
	V(x) = \frac{V_0}{2}\left[1+\tanh\left(\frac{x}{L}\right)\right], \quad x\in\mathbb{R},
\end{equation}
where $L$ controls the gradient and width of the step. The potential is continuous in order to avoid possible problems caused by discontinuity. Figure \ref{fig:KP_V} shows the electric potential $V(x)$ on $\Omega=(-20, 20)$, with $L=10^{-4}$ and $V_0 = 6.13\times 10^4$.

We take the initial condition
\begin{equation}\label{eq:KP_ini}
	\phi_1(0, x) = e^{ik_0x}e^{-\frac{(x-x_0)^2}{4}}, \quad 
	\phi_2(0, x) = Ce^{ik_0x}ee^{-\frac{(x-x_0)^2}{4}}, \quad x\in\mathbb{R},
\end{equation}
which represents the traveling Gaussian wave packet. The constant $C$ is given by
\begin{equation}
	C = \frac{ck_0}{mc^2+\sqrt{m^2c^4+c^2k_0^2}}.
\end{equation}
In the initial condition \eqref{eq:KP_ini}, $k_0$ stands for the wave packet momentum, and $x_0$ is the initial position. 

With this initial condition, the analytical transmission coefficient for the potential \eqref{eq:KP_potential} is \citep{KSG}
\begin{equation}
	{T}_\text{ana} = -\frac{\sinh(\pi kL)\sinh(\pi k'L)}{\sinh\left[\pi\left(\frac{V_0}{c}+k+k'\right)\frac{L}{2}\right]\sinh\left[\pi\left(\frac{V_0}{c}-k-k'\right)\frac{L}{2}\right]}, \quad V_0 > E_k + mc^2,
\end{equation}
where 
\begin{equation}
	k = \frac{1}{c}\sqrt{(E_k-V_0)^2 - m^2c^4}, \quad k' = -\frac{1}{c}\sqrt{E_k^2-m^2c^4},
\end{equation}
with
\begin{equation}
	E_k = \sqrt{k_0^2c^2 + m^2c^4}.
\end{equation}
In the computation, we take $k_0 = 106$, $L = 10^{-4}$, $x_0 = -10$. 

The simulation is computed until $T_\text{max} = 0.22$ on a bounded domain $x\in\Omega=(a, b)$, and periodic boundary conditions are assumed, which assure that the truncation error from the whole space problem is small enough to neglect. Take a positive even number $M$, define $h = (b-a)/M$ as the mesh size, and take $\tau>0$ to be the time step.

Denote $\Phi_f$ to be the outcome of the wave solution at $T_\text{max}$, then the numerical transmission coefficient is computed from $\Phi_f$ by
\begin{equation}
	{T}_\text{num} = \frac{\Phi_f^*\left((M/2+1):M\right)\Phi_f\left((M/2+1):M\right)}{\Phi_f^*\Phi_f}.
\end{equation}

\begin{figure}[htbp]
	\centering
	\begin{subfigure}{0.45\textwidth}
		\subcaption{The differences with different $\tau$}
		\includegraphics[width=1\textwidth]{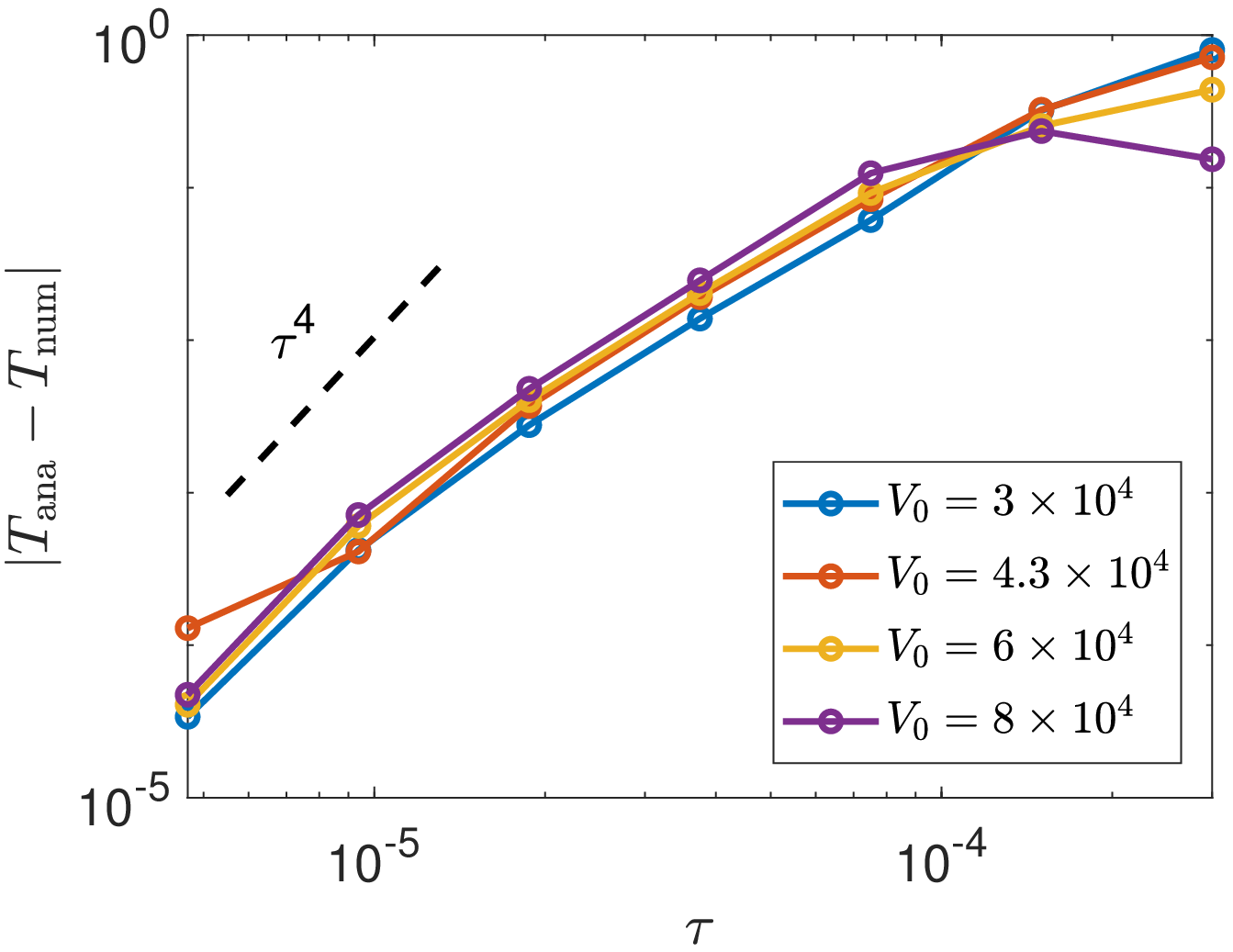}
	\end{subfigure}
    \begin{subfigure}{0.45\textwidth}
    	\subcaption{Comparison between $T_\text{ana}$ and $T_\text{num}$}
    	\includegraphics[width=1\textwidth]{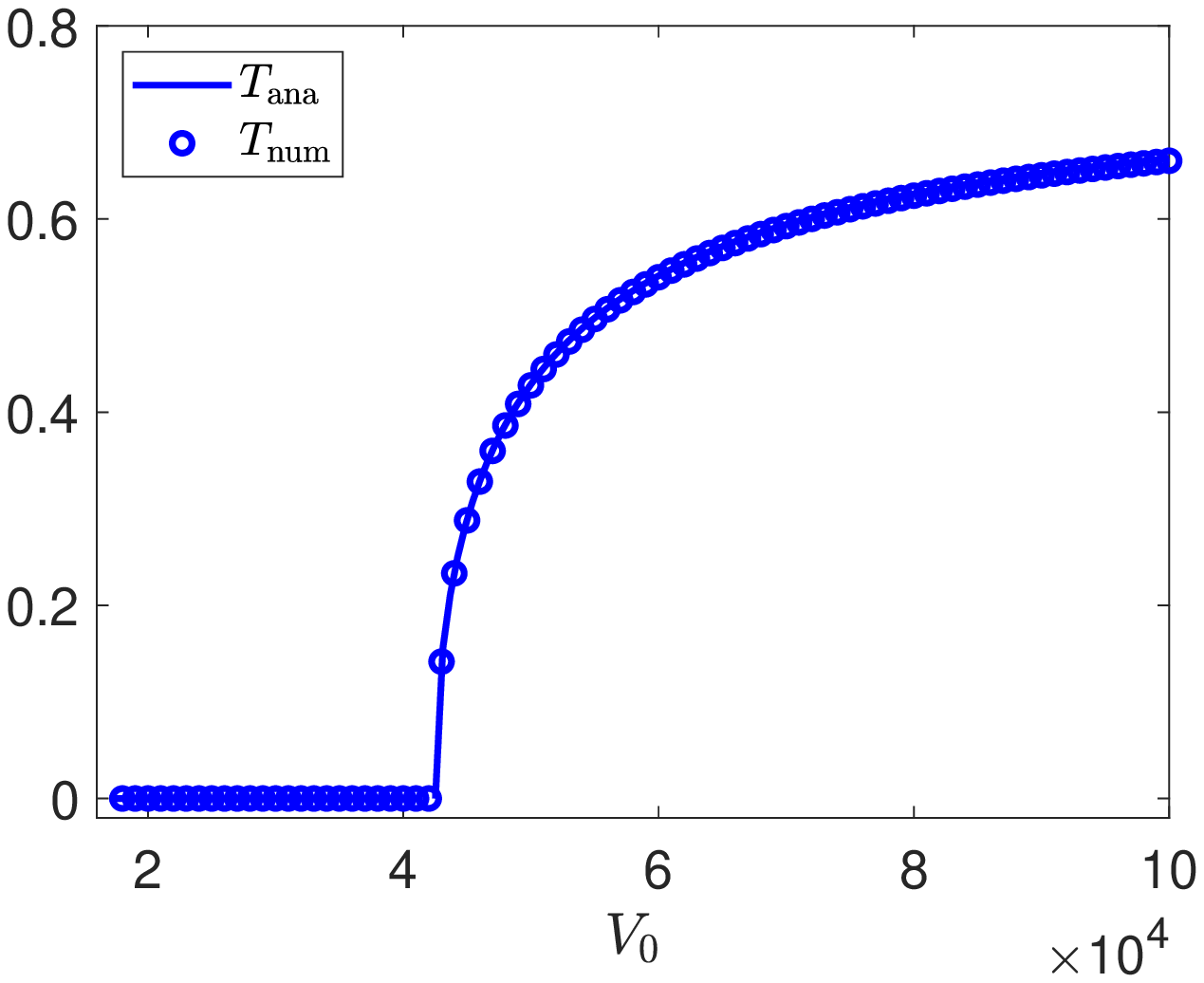}
    \end{subfigure}
    \caption{(a) The differences between $T_\text{ana}$ and $T_\text{num}$ from $S_\text{4c}$ with mesh size $h = 1/8192$ and different time step $\tau$; (b) Comparison between $T_\text{ana}$ and $T_\text{num}$ from $S_\text{4c}$ with mesh size $h = 1/2048$ and time step $\tau = 5\times 10^{-6}$.}
	\label{fig:KP_cmp}
\end{figure}

In this example, we take $a=-20$ and $b = 20$.
To show that $S_\text{4c}$ is fourth order accurate in time, we choose four different $V_0$, and fix the mesh size to be $h = 1/8192$. Differences between $T_\mathrm{ana}$ and $T_\mathrm{4c}$ are plotted in Figure \ref{fig:KP_cmp}(a), where we could observe that for small enough time step $\tau$, there is fourth order convergence. This validates that $S_\text{4c}$ is fourth order in time.

In addition, to verify the accuracy of $S_\text{4c}$, we compare the numerical results $T_\mathrm{num}$ with the analytical solution $T_\mathrm{ana}$ for different $V_0$. The mesh size here is fixed at $h = 1/2048$, which gives $M = 81920$ grid points, and the time step is taken to be $\tau = 5\times 10^{-6}$. Figure \ref{fig:KP_cmp}(b) exhibits the comparison between numerical transmission coefficients with analytical ones for different $V_0$.

The relative error of $T_\text{num}$ compared to $T_\text{ana}$ is always smaller than $0.4\%$ when $V_0>E_k+mc^2$. Additionally, when $V_0<E_k+mc^2$, $T_\text{num}$ is always nearly $0$, which corresponds well to the analytical analysis. These results suggest that our $S_\text{4c}$ scheme is accurate to solve the time-dependent Dirac equation.\\

In the following numerical examples, we consider the Dirac equation \eqref{eq:dirac2} with initial value \eqref{eq:initial11} on a bounded domain $\Omega$ with periodic boundary conditions.

We take mesh size $h>0$ in the numerical scheme, and apply Fourier spectral discretization in space, so that the steps involving $e^{\tau T}$ in \eqref{eq:S4c_td} could be easily solved in the phase space. The other steps involving $e^{\tau W}$ or $e^{\tau \widehat{W}}$ could be directly solved in the physical space.
Take time step size $\tau>0$ as before, then the temporal errors for the wave function, probability density and current density are respectively introduced as
\begin{align}
e_\Phi(t_n) = \left\|\Phi^n - \Phi(t_n, \cdot)\right\|_{l^2}, \; e_\rho(t_n) = \left\||\Phi^n|^2 - |\Phi(t_n, \cdot)|^2\right\|_{l^2}, \; e_\bJ(t_n) = \left\|\bJ(\Phi^n) - \bJ(\Phi(t_n, \cdot))\right\|_{l^2}
\end{align}
to represent the results, where $\bJ(\Phi) = (\bJ_1(\Phi), \bJ_2(\Phi))^T$, and 
\be
\bJ_l(\Phi) = (\Phi)^*\sigma_l\Phi,  \quad l = 1, 2.
\ee

\subsection{An example in 1D}

In the example, we take $d = 1$ in \eqref{eq:dirac2}, and the initial conditions are set to be
\be
\phi_1(0, x) = e^{-x^2/2}, \quad \phi_2(0, x) = e^{-(x-1)^2/2}, \quad x\in\mathbb{R}.
\ee
The time-dependent electromagnetic potentials are taken as
\be
V(t, x) = \frac{1-tx}{1+t^2x^2}, \quad A_1(t,x) = \frac{(tx+1)^2}{1+t^2x^2}, \quad t>0, \quad x\in\mathbb{R}.
\ee
The problem is solved numerically on a bounded domain $\Omega = (-32, 32)$. As the analytical solution is unavailable, to obtain the `exact' solution, fine mesh size $h_e = 1/16$ and fine time step size $\tau_e = 10^{-5}$ are used in $S_\text{4c}$ \eqref{eq:S4c_td}. 

The temporal errors in this example are quantified as
\begin{align*}
&e_\Phi(t_n) = \left\|\Phi^n - \Phi(t_n, \cdot)\right\|_{l^2} := \sqrt{h\sum_{j=0}^{M-1}|\Phi_j^n - \Phi(t_n, x_j)|^2},\\
&e_\rho(t_n) =  \left\||\Phi^n|^2 - |\Phi(t_n, \cdot)|^2\right\|_{l^2} := \sqrt{h\sum_{j=0}^{M-1}\left(|\Phi_j^n|^2 - |\Phi(t_n, x_j)|^2\right)^2},\\
&e_\bJ(t_n) = \left\|\bJ(\Phi^n) - \bJ(\Phi(t_n, \cdot))\right\|_{l^2} := \sqrt{h\sum_{j=0}^{M-1}\sum_{k=1}^2\left|(\Phi_j^n)^*\sigma_k\Phi_j^n - (\Phi(t_n, x_j))^*\sigma_k\Phi(t_n, x_j)\right|^2},
\end{align*}
with $M = 64/h$, $x_j := -32+jh$, $j = 0, ..., M$, and the numerical solution $\Phi^n := (\Phi_0^n, \Phi_1^n, ..., \Phi_{M-1}^n)^T$.

\begin{figure}[htp!]
	\centering
	\begin{subfigure}{0.3\textwidth}
		\caption{Errors for the wave fucntion}
	    \includegraphics[width=\textwidth]{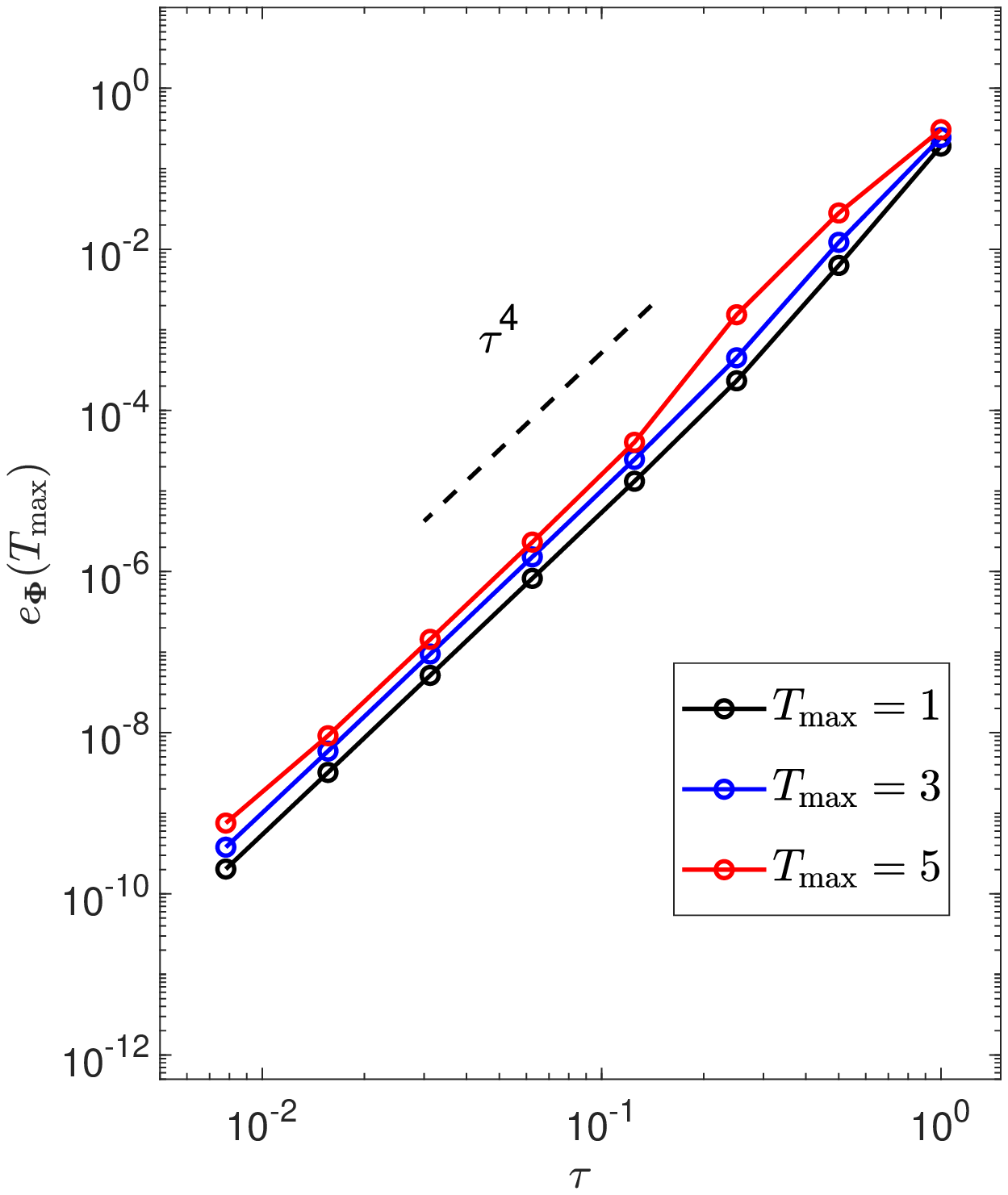}
	\end{subfigure}
    \begin{subfigure}{0.3\textwidth}
    	\caption{Errors for the probability density}
	    \includegraphics[width=\textwidth]{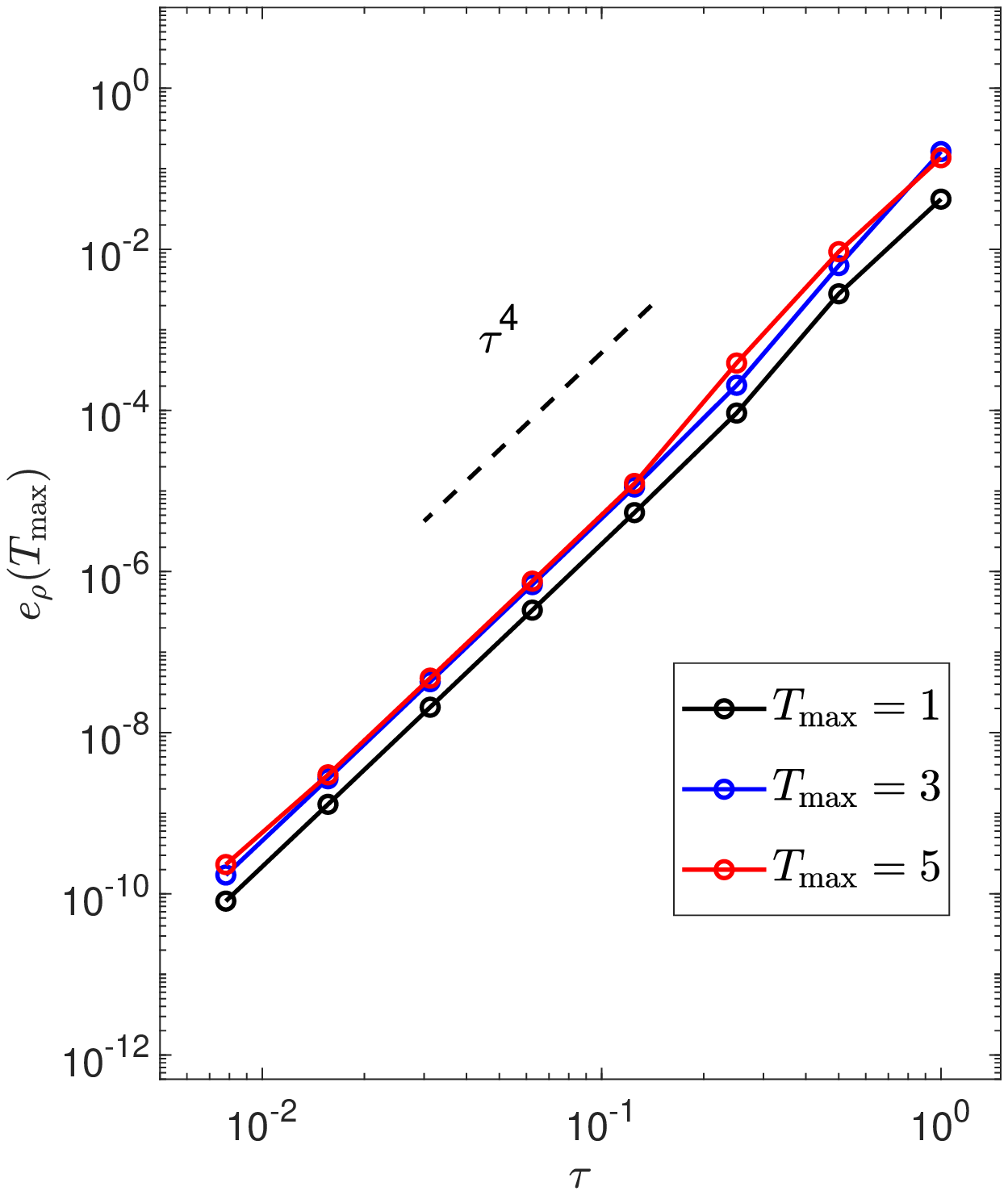}
	\end{subfigure}
    \begin{subfigure}{0.3\textwidth}
    	\caption{Errors for the current density}
    	\includegraphics[width=\textwidth]{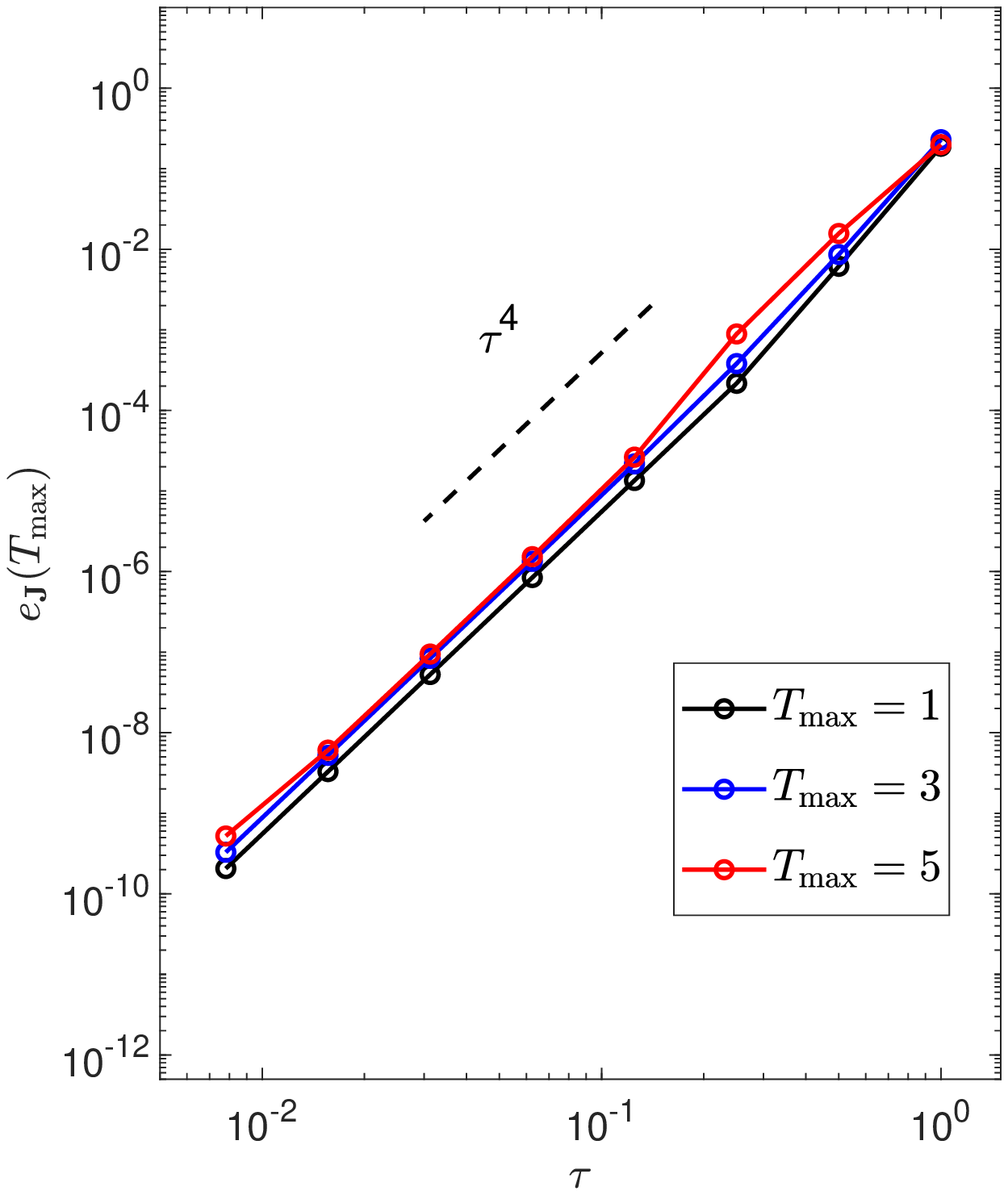}
	\end{subfigure}
	\caption{Temporal errors for the wave function, probability density, and current density with different $T_\mathrm{max}$s, 1D case.}
	\label{fig:1D_err}
\end{figure}

Figure \ref{fig:1D_err} shows $e_\Phi(T)$, $e_\rho(T)$ and $e_\bJ(T)$ respectively for different $T_\mathrm{max}$s.

From the figure, we could clearly observe fourth order convergence in time for the wave function, probability density and current density by applying $S_\text{4c}$ \eqref{eq:S4c_td} to the Dirac equation in 1D with time-dependent potentials. When $T_\mathrm{max}$ becomes larger, there is a slight increase in the error for a fixed time step size, and the performance for large time step sizes is influenced by a bit. But overall, the fourth-order convergence is not affected. Consequently, $S_\text{4c}$ \eqref{eq:S4c_td} performs well in this 1D case.\\

Additionally, in order to compare the performance of different splitting methods, we also apply the first-order ($S_1$) \citep{Trotter}, the second-order ($S_2$) \citep{Strang}, the fourth-order Forest-Ruth  ($S_4$) \citep{FR, Suzuki0, Yoshida}, and the fourth-order Runge-Kutta $S_\text{4RK}$ \citep{Geng} splitting methods to the Dirac equation with time-dependent potentials. The ideas of application are similar to $S_\text{4c}$, where we use the time-ordering technique. To observe the results more clearly, we take the bounded domain $\Omega = (-64, 64)$, and the fine mesh size $h_e = 1/ 64$. The initial value and electromagnetic potentials are taken as before.

The results from the five splitting methods are summarized in Table \ref{Table 1}.

 \begin{table}[h!]\renewcommand{\arraystretch}{1.6}
	\small
	\centering
	\begin{tabular}{|c|c|ccccccc|}
		\hline
		\multicolumn{2}{|c|}{}&  $\tau_0 = 1 / 2$ & $\tau_0 / 2$ & $\tau_0 / 2^2$ &
		$\tau_0 / 2^3$ & $\tau_0 / 2^4$ & $\tau_0 / 2^5$  & $\tau_0 / 2^6$ \\
		\hline
		\multirow{3}{*}{$S_1$}&  $e_\Phi(t=5)$ &  9.25E-1 &	3.60E-1 & 1.61E-1 &	7.72E-2 & 3.79E-2 & 1.88E-2 & 9.37E-3 \\
		%\cline{3-10}
		&rate  & --	 & 	1.36 &	1.16 &	1.06 &	1.03 &	1.01 &	1.01 \\
		%\cline{3-10}
		&CPU Time  &  0.03 &	0.05 &	0.10 &	0.13 &	0.19 &	0.40 &	\textbf{0.75} \\
		\hline
		\multirow{3}{*}{$S_2$} & $e_\Phi(t=5)$  &  	6.14E-1 &	1.51E-1 &	3.76E-2 &	9.39E-3 &	2.35E-3 &	5.87E-4 &	1.47E-4 \\
		%\cline{2-10}
		& rate  & -- & 	2.03 &	2.00 &	2.00 &	2.00 &	2.00 &	2.00 \\
		%\cline{2-10}
		& CPU Time  &  0.05 &	0.06 &	0.12 &	0.13 &	0.25 &	0.53  &	\textbf{1.10} \\
		\hline
		\multirow{3}{*}{$S_4$}&$e_\Phi(t=5)$  &   2.21E-1 &	2.37E-2 &	1.82E-3 &	1.22E-4 &	7.80E-6 &	4.90E-7 &	3.07E-8 \\
		%\cline{2-10}
		&rate  & -- &	3.22 &	3.70 &	3.89 &	3.97 &	3.99 &	4.00 \\
		%\cline{2-10}
		&CPU Time  & 	0.10 &	0.12 &	0.22 &	0.38 &	0.78 &	1.38 &	\textbf{2.89} \\
		\hline
		\multirow{3}{*}{$S_\text{4c}$} & $e_\Phi(t=5)$  &  2.82E-2 &	1.54E-3 &	4.04E-5 &	2.32E-6 &	1.44E-7 &	8.95E-9 &	5.94E-10 \\
		%\cline{2-10}
		& rate  & -- & 	4.19 &	5.26 &	4.12 &	4.02 &	4.00 &	3.91 \\
		%\cline{2-10}
		& CPU Time  &  	0.07 &	0.09 &	0.13 &	0.25 &	0.45 &	0.88 &	\textbf{1.78} \\
		\hline
		\multirow{3}{*}{$S_\text{4RK}$} & $e_\Phi(t=5)$  &   4.25E-3 &	2.11E-4 &	7.42E-6 &	4.52E-7 &	2.82E-8 &	1.78E-9 &	2.15E-10 \\
		%\cline{2-10}
		& rate  & -- &  4.33 &	4.83 &	4.04 &	4.00 &	3.99 &	3.05 \\
		%\cline{2-10}
		& CPU Time  &  0.11 &	0.16 &	0.29 &	0.59 &	1.10 &	2.27 &	\textbf{5.32} \\
		\hline
	\end{tabular}
    \vspace{10pt}
	\caption{Temporal errors $e_\Phi(t=5)$ of different time-splitting methods under different time step sizes $\tau$ for the  Dirac equation \eqref{eq:dirac2} in 1D. Here we also list convergence rates and computational time  (CPU time in seconds) for comparison.}
	\label{Table 1}
\end{table}

Because the convergence behaviors of the errors for wave function, probability density and current density are similar, here we only list the results for $e_\Phi(t)$. From Table \ref{Table 1}, we can see that these methods all achieve expected order of convergence. Similar to the case with time-independent electromagnetic potentials, the computational costs for the three fourth-order methods $S_4$, $S_\text{4c}$, $S_\text{4RK}$ are approximately three times, twice, and five to six times the time costs for $S_1$ and $S_2$, respectively. In this sense, $S_\text{4c}$ performs much better than the other two methods. Moreover, under the same time step size, the error $e_\Phi(t=5)$ for $S_\text{4c}$ is comparable to the error for $S_\text{4RK}$, and is about 50 times smaller than the error for $S_4$. Consequently, we conclude that $S_\text{4c}$ is efficient and accurate for the Dirac equation with time-dependent potentials, and is the best to apply among the three fourth-order methods.\\

To show that $S_\text{4c}$ \eqref{eq:S4c_td} is still valid for higher dimensions, we give examples in 2D as follows.

\subsection{Examples in 2D}
In the 2D examples, we take $d=2$ in \eqref{eq:dirac2}, and give the initial data:
\begin{equation}
\phi_1(0, \bx) =  e^{-\frac{x^2 + y^2}{2}}, \quad \phi_2(0, \bx) = e^{-\frac{(x - 1)^2 + y^2}{2}}, \qquad \bx=(x, y)^T\in\mathbb{R}^2.
\end{equation}
The time-dependent potentials are taken in honey-comb form
\be\begin{split}
	\label{potential_2D}
	V(t, \mathbf{x}) &= \cos\left(\frac{4\pi}{\sqrt{3}}\mathbf{e_1}(t)\cdot \mathbf{x}\right) +
	\cos\left(\frac{4\pi}{\sqrt{3}}\mathbf{e}_2(t)\cdot\mathbf{x}\right) + \cos\left(\frac{4\pi}{\sqrt{3}}\mathbf{e}_3(t)
	\cdot \mathbf{x}\right),\\
	A_1(t, \mathbf{x}) &= A_2(t, \mathbf{x}) = 0,\qquad  \bx\in{\mathbb R}^2,
\end{split}
\ee
with
\begin{equation}
\begin{aligned}
&\mathbf{e}_1(t) = (\cos(\theta(t)), \sin(\theta(t)))^T, \quad \mathbf{e}_2(t) = (\cos(\theta(t)+\frac{2\pi}{3}), \sin(\theta(t) + \frac{2\pi}{3}))^T,\\
&\qquad\qquad \mathbf{e}_3(t) = (\cos(\theta(t)+\frac{4\pi}{3}), \sin(\theta(t)+\frac{4\pi}{3}))^T,
\end{aligned}
\end{equation}
where $\theta(t)$ is a given function. In our examples, we consider $\theta(t)$ to be
\begin{flalign*}
& \quad\text{(1)} \quad\theta(t) \equiv \pi; &\\
& \quad\text{(2)}\quad\theta(t) = \pi + \pi t; &\\
& \quad\text{(3)}\quad\theta(t) = \pi + \pi\cos(\pi t). &
\end{flalign*}

\begin{figure}[!htbp]
	\centering
	\includegraphics[width=0.9\textwidth]{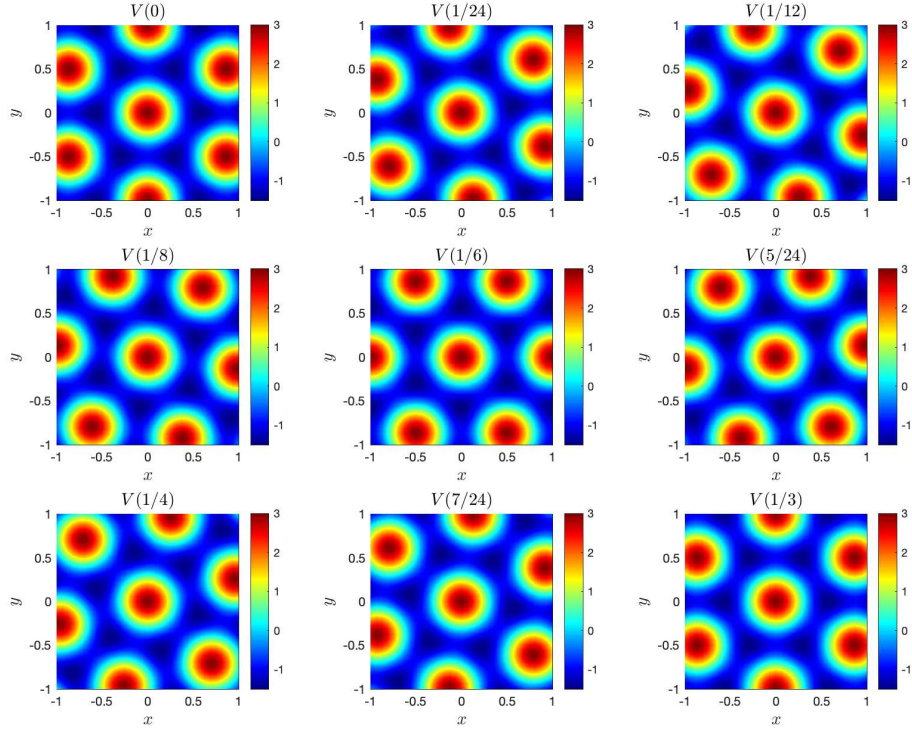}
	\caption{The potential $V(t)$ with $\theta(t) = \pi + \pi t$ from $t=0$ to $t=1/3$.}
	\label{fig:V_ic2}
\end{figure}

\begin{figure}[!htbp]
	\centering
	\includegraphics[width=0.9\textwidth]{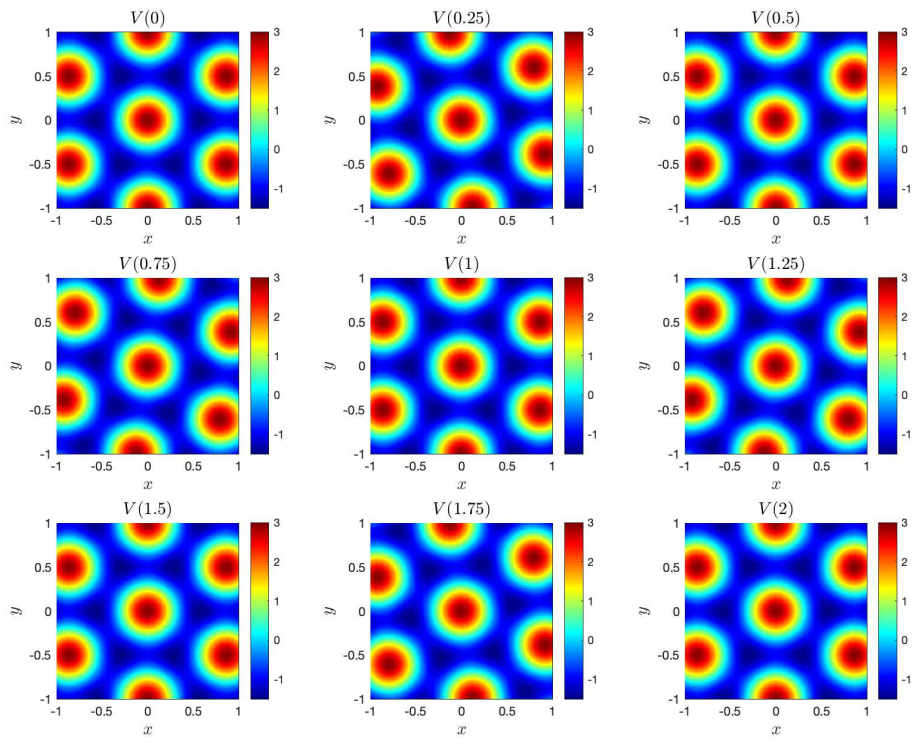}
	\caption{The potential $V(t)$ with $\theta(t) = \pi + \pi\cos(\pi t)$ from $t=0$ to $t=2$.}
	\label{fig:V_ic3}
\end{figure}

The varying potentials in cases (2) and (3) are illustrated in Figure \ref{fig:V_ic2} and \ref{fig:V_ic3}, respectively. Here we take $V(t) := V(t, \cdot)$ for short. As the potentials are periodic in space, only those in domain $[-1, 1]\times[-1, 1]$ are exhibited for better illustration. The potential in case (1) is fixed as $V(0)$ in case (2) (cf. Figure \ref{fig:V_ic2}).

Through simple computation, we could get the period in time of case (2) is $1/3$, and the period in time of case (3) is $2$, which corresponds well with the figures. Indeed, in case (2), there is anticlockwise rotation of the local circle potentials with respect to the center $(0, 0)$, and after $\triangle t = 1/3$, the circle potentials are all back to the initial positions. In case (3), the local circle potentials would oscillate along a circle centered at $(0, 0)$, and $t=2$ is when the first period ends.

We set the magnetic potentials to $0$ so that $S_\text{4c}$ \eqref{eq:S4c_td} could be efficiently applied. The problem is solved numerically on a bounded domain $\Omega = (-25, 25)\times (-25, 25)$.

Similar to the 1D example, we obtain a numerical `exact' solution by using the $S_\text{4c}$ \eqref{eq:S4c_td} with a
fine mesh size $h_e = \frac{1}{16}$ and a small time step $\tau_e = 10^{-4}$.

The temporal errors in this example are quantified as
\begin{align*}
&e_\Phi(t_n) = \left\|\Phi^n - \Phi(t_n, \cdot)\right\|_{l^2} := {h\sqrt{\sum_{j = 0}^{M - 1}\sum_{l = 0}^{M - 1}|\Phi_{jl}^n -  \Phi(t_n, x_j,y_l)|^2}},\\
&e_\rho(t_n) =  \left\||\Phi^n|^2 - |\Phi(t_n, \cdot)|^2\right\|_{l^2} := {h\sqrt{\sum_{j=0}^{N-1}\sum_{l = 0}^{M - 1}\left(|\Phi_{jl}^n|^2 - |\Phi(t_n, x_j, y_l)|^2\right)^2}},\\
&e_\bJ(t_n) = \left\|\bJ(\Phi^n) - \bJ(\Phi(t_n, \cdot))\right\|_{l^2} := {h\sqrt{\sum_{j=0}^{N-1}\sum_{l = 0}^{M - 1}\sum_{k=1}^2\left|(\Phi_{jl}^n)^*\sigma_k\Phi_{jl}^n - (\Phi(t_n, x_j, y_l))^*\sigma_k\Phi(t_n, x_j, y_l)\right|^2}},
\end{align*}
with $M = 50/h$, $x_j := -25+jh$, $y_l := -25+lh$, and $\Phi_{jl}^n$ is the numerical solution at $(x_j, y_l)$ for time $t = n\tau$. Here $j$, $l = 0, ..., M$, $n = 0, 1, ..., T/\tau$. We show the results case by case.\\

(1) $\theta(t)\equiv \pi$.

In this case, $\theta(t)$ is time-independent, so that the method is equivalent to $S_\text{4c}$ for the Dirac equation with time-independent potentials \citep{BY}. The results for $e_\Phi(t=3)$, $e_\rho(t=3)$, and $e_\bJ(t=3)$ are shown in Table \ref{tb:2D_ic1}.

\begin{table}[h!]\renewcommand{\arraystretch}{1.6}
	\small
	\centering
	\begin{tabular}{|c|cccccccc|}
		\hline
		& $\tau_0 = 1 / 2$ & $\tau_0 / 2$ & $\tau_0 / 2^2$ &
		$\tau_0 / 2^3$ & $\tau_0 / 2^4$ & $\tau_0 / 2^5$  & $\tau_0 / 2^6$ & $\tau_0 / 2^7$  \\
		\hline
	    $e_\Phi(t=3)$ &  2.13E-1 &	9.67E-3 &	2.37E-4 &	1.41E-5 &	8.76E-7 &	5.46E-8  &	3.41E-9 &	2.14E-10 \\
		rate  & --	 & 	4.46 &	5.35 &	4.07 &	4.01 &	4.00 &	4.00 &	3.99\\
		\hline
	    $e_\rho(t=3)$ & 1.04E-1 &	3.86E-3 &	7.91E-5 &	4.63E-6 &	2.86E-7 &	1.78E-8 &	1.11E-9 &	7.02E-11\\
		rate  & -- &  4.75 &	5.61 &	4.10 &	4.02 &	4.00 &	4.00 &	3.98\\
		\hline
		$e_\bJ(t=3)$ & 1.28E-1 &	5.60E-3 &	1.13E-4 &	6.70E-6 &	4.15E-7 &	2.59E-8 &	1.62E-9 &	1.04E-10\\
		rate &  -- &	4.51 &	5.63 &	4.07 &	4.01 &	4.00 &	4.00 &	3.96\\
		\hline
	\end{tabular}
	\vspace{10pt}
	\caption{Temporal errors $e_\Phi(t=3)$, $e_\rho(t=3)$, and $e_\bJ(t=3)$ for the  Dirac equation \eqref{eq:dirac2} in 2D, with the potential given in \eqref{potential_2D}, where $\theta(t) \equiv \pi$.}
	\label{tb:2D_ic1}
\end{table}

\begin{figure}[!htbp]
	\centering
	\includegraphics[width=0.9\textwidth]{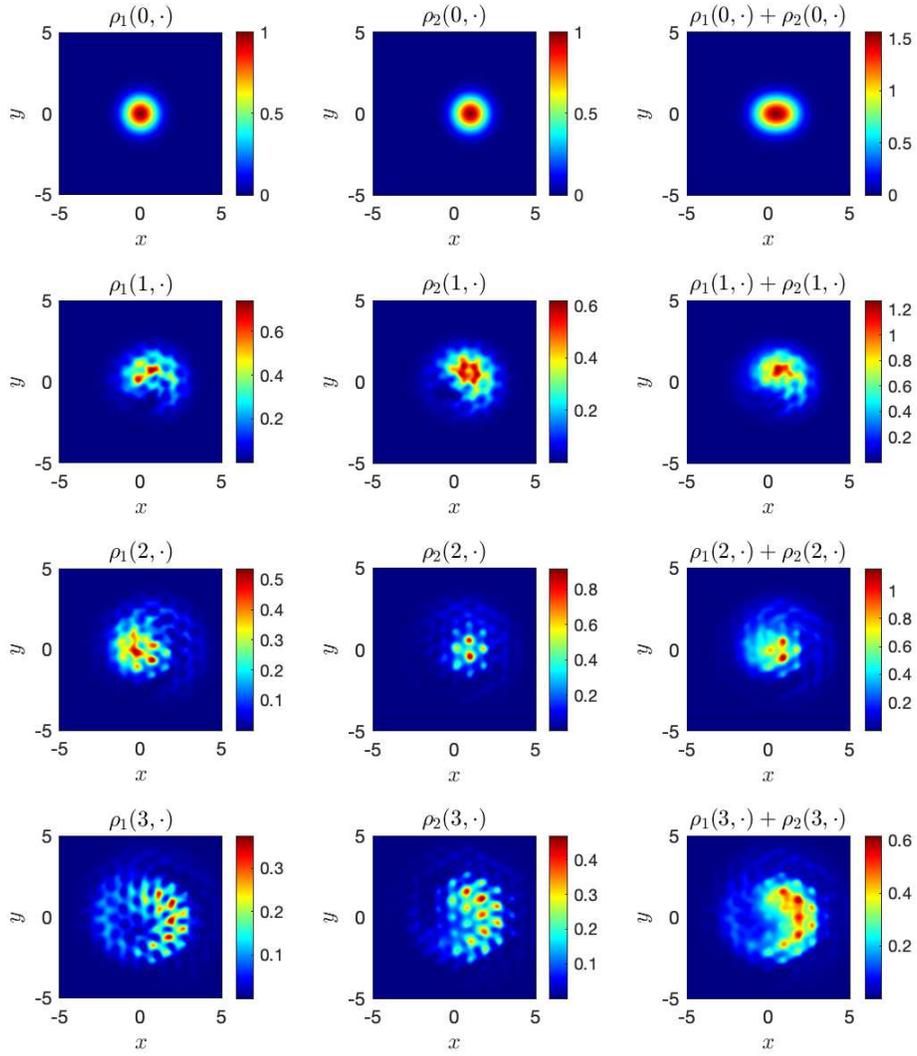}
	\caption{The probability densities $\rho_1(t, \cdot)$, $\rho_2(t, \cdot)$, and their sum $\rho_1(t, \cdot) + \rho_2(t, \cdot)$ with $t = 0$, $1$, $2$, $3$, when $\theta(t)\equiv\pi$.}
	\label{fig:ic1}
\end{figure}

From the table, we could observe clear fourth-order convergence for the wave function, probability density, and current density. The evolution of $\rho_1(t) := \rho_1(t, \bx)$, $\rho_2(t) := \rho_2(t, \bx)$, which respectively represents the probability density of the two components, and their sum is shown in Figure \ref{fig:ic1}.\\

(2) $\theta(t) = \pi + \pi t$.

In this case, $\theta(t)$ is monotonically increasing, which results in a periodic electric potential $V(t) := V(t, \bx)$. Table \ref{tb:2D_ic2} gives $e_\Phi(t=3)$, $e_\rho(t=3)$, and $e_\bJ(t=3)$ under this potential.

\begin{table}[h!]\renewcommand{\arraystretch}{1.6}
	\small
	\centering
	\begin{tabular}{|c|ccccccc|}
		\hline
		& $\tau_0 = 1 / 8$ & $\tau_0 / 2$ & $\tau_0 / 2^2$ &
		$\tau_0 / 2^3$ & $\tau_0 / 2^4$ & $\tau_0 / 2^5$  & $\tau_0 / 2^6$  \\
		\hline
		$e_\Phi(t=3)$ &  5.09E-1 &	6.61E-2 &	2.69E-4 &	1.31E-5 &	7.79E-7 &	4.81E-8 &	3.00E-9 \\
		rate  & --	 & 	2.95 &	7.94 &	4.36 &	4.07 &	4.02 &	4.00\\
		\hline
		$e_\rho(t=3)$ & 1.07E-1 &	3.51E-3 &	1.11E-5 &	6.49E-7 &	4.00E-8 &	2.49E-9 &	1.56E-10\\
		rate  & -- &  4.93 &	8.31 &	4.09 &	4.02 &	4.01 &	4.00\\
		\hline
		$e_\bJ(t=3)$ & 1.54E-1 &	5.59E-3 &	1.82E-5 &	1.03E-6 &	6.29E-8 &	3.91E-9 &	2.44E-10\\
		rate &  -- &	4.79 &	8.26 &	4.15 &	4.03 &	4.01 &	4.00\\
		\hline
	\end{tabular}
	\vspace{10pt}
	\caption{Temporal errors $e_\Phi(t=3)$, $e_\rho(t=3)$, and $e_\bJ(t=3)$ for the  Dirac equation \eqref{eq:dirac2} in 2D, with the potential given in \eqref{potential_2D}, where $\theta(t) = \pi + \pi t$.}
	\label{tb:2D_ic2}
\end{table}

\begin{figure}[!htbp]
	\centering
	\includegraphics[width=0.9\textwidth]{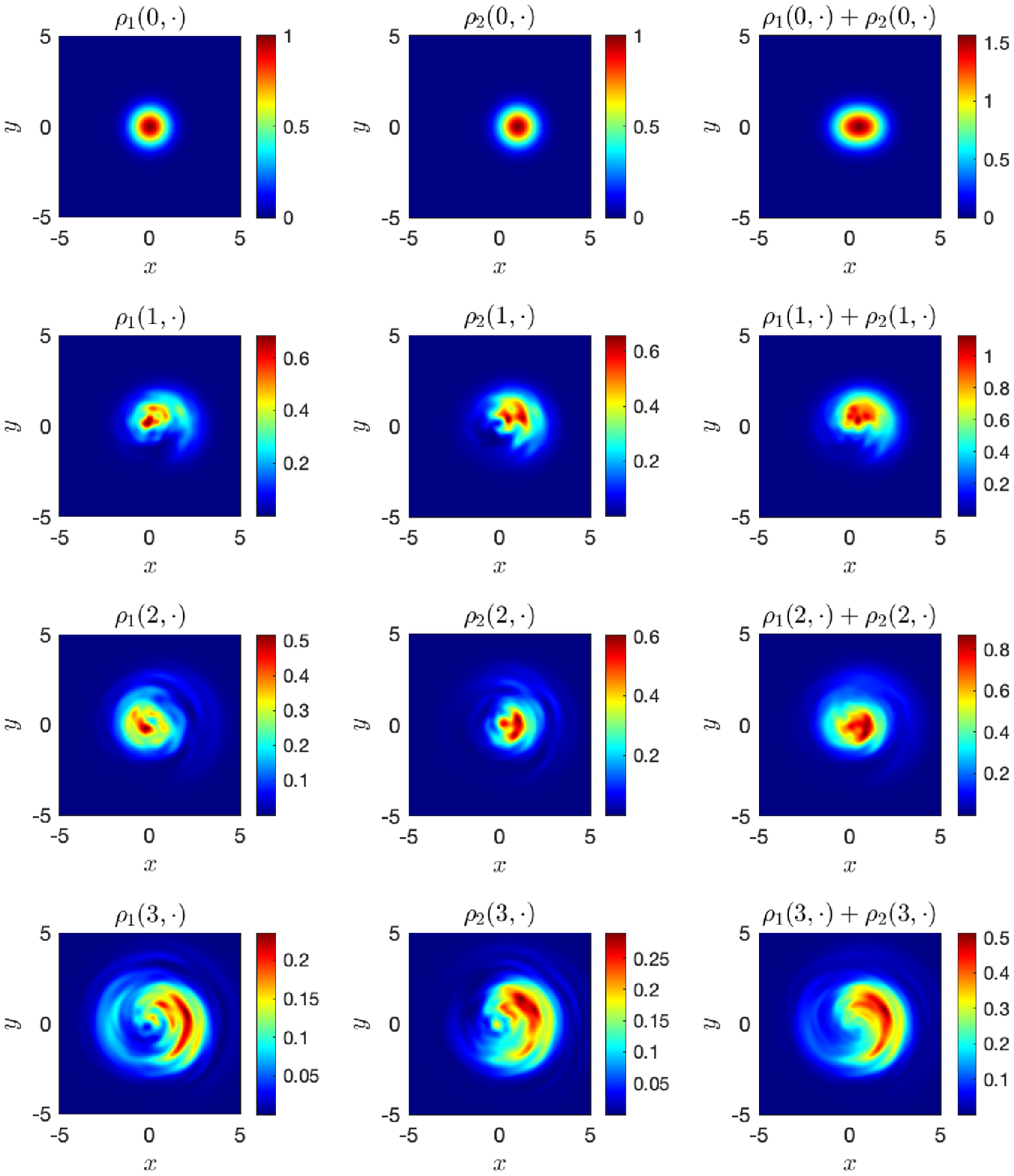}
	\caption{The probability densities $\rho_1(t, \cdot)$, $\rho_2(t, \cdot)$, and their sum $\rho_1(t, \cdot) + \rho_2(t, \cdot)$ with $t = 0$, $1$, $2$, $3$, when $\theta(t) = \pi + \pi t$.}
	\label{fig:ic2}
\end{figure}

From the table, we could observe that when the time step size is large, there is no fourth-order convergence. But by further decreasing time step sizes, we would obtain fourth-order convergence for the wave function and the two physical observables, which validates $S_\text{4c}$ \eqref{eq:S4c_td} with time-dependent potential for the Dirac equation in 2D. The dynamics of $\rho_1(t)$, $\rho_2(t)$, and their sum in this case is given in Figure \ref{fig:ic2}.\\

(3) $\theta(t) = \pi + \pi\cos(\pi t)$.

In this case, $\theta(t)$ is periodic in time, which generates a periodic electric potential $V(t)$ with the same period. Table \ref{tb:2D_ic3} gives $e_\Phi(t=3)$, $e_\rho(t=3)$, and $e_\bJ(t=3)$ under this potential.

\begin{table}[h!]\renewcommand{\arraystretch}{1.6}
	\small
	\centering
	\begin{tabular}{|c|ccccccc|}
		\hline
		& $\tau_0 = 1 / 8$ & $\tau_0 / 2$ & $\tau_0 / 2^2$ &
		$\tau_0 / 2^3$ & $\tau_0 / 2^4$ & $\tau_0 / 2^5$  & $\tau_0 / 2^6$  \\
		\hline
		$e_\Phi(t=3)$ &  8.53E-1 &	2.74E-1 &	4.08E-2 &	2.48E-3 &	3.92E-8 &	2.45E-9 &	1.54E-10 \\
		rate  & --	 & 1.64	& 2.75 &	4.04 &	15.95 &	4.00 &	3.99\\
		\hline
		$e_\rho(t=3)$ & 2.46E-1 &	7.35E-2 &	7.65E-3 &	5.77E-5 &	6.51E-9 &	4.05E-10 &	2.60E-11\\
		rate  & -- &  1.74 &	3.26 &	7.05 &	13.11 &	4.01 &	3.96\\
		\hline
		$e_\bJ(t=3)$ & 3.68E-1 &	1.07E-1 &	1.07E-2 &	9.03E-5 &	1.18E-8 &	7.28E-10 &	4.54E-11\\
		rate &  -- & 1.78 &	3.33 &	6.88 &	12.90 &	4.02 &	4.00\\
		\hline
	\end{tabular}
	\vspace{10pt}
	\caption{Temporal errors $e_\Phi(t=3)$, $e_\rho(t=3)$, and $e_\bJ(t=3)$ for the  Dirac equation \eqref{eq:dirac2} in 2D, with the potential given in \eqref{potential_2D}, where $\theta(t) = \pi + \pi\cos(
		\pi t)$.}
	\label{tb:2D_ic3}
\end{table}

\begin{figure}[!htbp]
	\centering
	\includegraphics[width=0.9\textwidth]{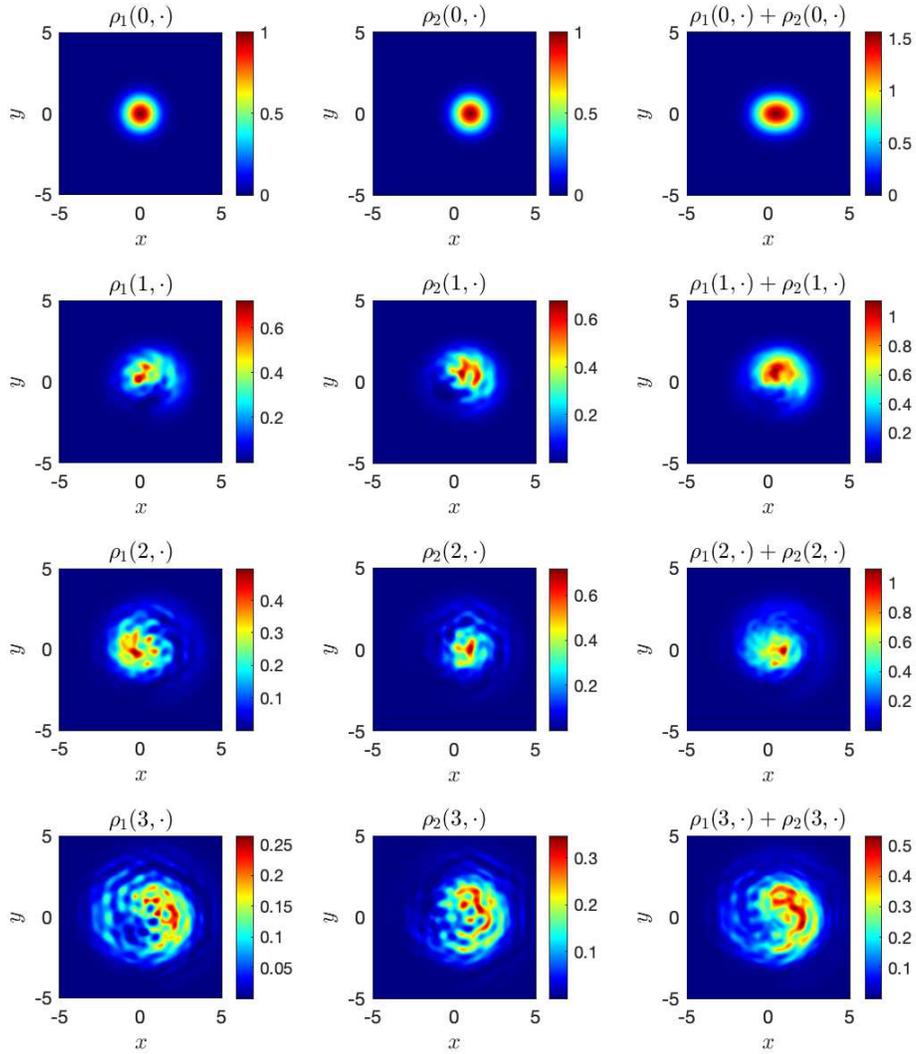}
	\caption{$\rho_1(t)$, $\rho_2(t)$, and $\rho_1(t) + \rho_2(t)$ with $t = 0$, $1$, $2$, $3$, where $\theta(t) = \pi + \pi\cos(\pi t)$.}
	\label{fig:ic3}
\end{figure}

The conclusions we could draw from this table is similar to case (2). When the time step size is large, the fourth-order convergence is not obtained. When the time step size is small enough, we could observe fourth-order convergence, which again validates $S_\text{4c}$ \eqref{eq:S4c_td} for time-dependent potentials. The dynamics of $\rho_1(t)$, $\rho_2(t)$, and their sum in this case is given in Figure \ref{fig:ic3}.\\

Overall, from the three numerical examples, we could conclude that the $S_\text{4c}$ derived for the Dirac equation with time-dependent potentials is valid in 2D. It is simple to apply when there is no magnetic potentials, and the results are satisfactory. The method successfully captures different dynamics of the probability densities under various electric potentials.\\

\section{Conclusion}
In this paper, we study the fourth-order compact time-splitting method ($S_\text{4c}$) for the Dirac equation with time-dependent potentials. The time-ordering technique is introduced to deal with the time-dependence, so that in each time step, the choices of $t$ for those sub-steps with potentials vary. Under this treatment, $S_\text{4c}$ remains efficient, as the overall computational cost does not increase much compared to the case with time-independent potentials. Numerical examples in 1D and 2D are given to validate the accuracy, and comparison of $S_\text{4c}$ with other splitting methods $S_1$, $S_2$, $S_4$, $S_\text{4RK}$ is also exhibited, which shows that $S_\text{4c}$ performs the best considering efficiency and accuracy.\\

\noindent{\bf Acknowledgments} This work was partially supported by the Ministry of Education of Singapore grant R-146-000-247-114. Part of the work was done when the author was visiting the Institute for Mathematical Sciences at the National University of Singapore in 2020. The author is grateful to Prof. Weizhu Bao at National University of Singapore for fruitful discussions.\\

\newpage
{\bf Appendix A}. Derivation of the double commutator in Lemma \ref{lemma:db_com_1d} for the Dirac equation \eqref{eq:dirac2} in 1D.\\
\setcounter{equation}{0}
\renewcommand{\theequation}{A.\arabic{equation}}
\indent It is easy to check that $[W(t), [T_1+T_2, W(t)]] = [W(t), [T_1, W(t)]] + [W(t), [T_2, W(t)]]$. Based on this relation, the double commutators in 1D can be derived as follows.\\

From \eqref{eq:split_TW}, in 1D, we have
\be
T = -\sigma_1\partial_1 - i\sigma_3, \quad W(t) = -i\left(V(t, x)I_2 - A_1(t, x)\sigma_1\right).
\ee
Through the linearity of the double commutator in $T$,
\begin{equation}\label{eq:TWd1}
	[W(t), [T, W(t)]] = -{\left[W(t), \left[\sigma_1\partial_1, W(t)\right]\right]} -
	i{\left[W(t), \left[\sigma_3, W(t)\right]\right]}.
\end{equation}
The two terms on the right hand side give
	\bea
	\label{Ts1d1}
	{\left[W(t), \left[\sigma_1\partial_1, W(t)\right]\right]}&=&2\left(-i\big(V(t, x)I_2 - A_1(t, x)\sigma_1\big)\right)
	\left(\sigma_1\partial_1\right)\left(-i\big(V(t, x)I_2 - A_1(t, x)\sigma_1\big)\right) \nn \\
	&& - \left(-i\big(V(t, x)I_2 - A_1(t, x)\sigma_1\big)\right)^2\left(\sigma_1\partial_1\right)
	- \left(\sigma_1\partial_1\right)\left(-i\big(V(t, x)I_2 - A_1(t, x)\sigma_1\big)\right)^2\nn\\
	&=& -2\big(V(t, x)I_2 - A_1(t, x)\sigma_1\big)\sigma_1\partial_1\big(V(t, x)I_2 - A_1(t, x)\sigma_1\big)\nn\\
	&& +\big(V(t, x)I_2 - A_1(t, x)\sigma_1\big)^2\sigma_1\partial_1
	+ \sigma_1\partial_1\big(V(t, x)I_2 - A_1(t, x)\sigma_1\big)^2\nn\\
	&=& -2\sigma_1\big(V(t, x)I_2 - A_1(t, x)\sigma_1\big)\partial_1\big(V(t, x)I_2 - A_1(t, x)\sigma_1\big)\nn\\
	&& -2\sigma_1\big(V(t, x)I_2 - A_1(t, x)\sigma_1\big)^2\partial_1 +2\sigma_1
	\big(V(t, x)I_2 - A_1(t, x)\sigma_1\big)^2\partial_1\nn\\
	&& + 2\sigma_1\big(V(t, x)I_2 - A_1(t, x)\sigma_1\big)\partial_1\big(V(t, x)I_2 - A_1(t, x)\sigma_1\big)\nn\\
	&=& 0,
	\eea
	and
	\bea
	\label{Ts1d2}
	{[W(t), [\sigma_3, W(t)]]}
	&=& 2\left(-i\big(V(t, x)I_2 - A_1(t, x)\sigma_1\big)\right)\sigma_3\left(-i
	\big(V(t, x)I_2 - A_1(t, x)\sigma_1\big)\right)\nn \\ 
	&&-\left(-i\big(V(t, x)I_2 - A_1(t, x)\sigma_1\big)\right)^2\sigma_3
	- \sigma_3\left(-i\big(V(t, x)I_2 - A_1(t, x)\sigma_1\big)\right)^2\nn\\
	&=& -2\big(V(t, x)I_2 - A_1(t, x)\sigma_1\big)\big(V(t, x)I_2 + A_1(t, x)\sigma_1\big)\sigma_3
	+ \big(V(t, x)I_2 - A_1(t, x)\sigma_1\big)^2\sigma_3\nn\\
	&& + \big(V(t, x)I_2 + A_1(t, x)\sigma_1\big)^2\sigma_3\nn\\
	&=& -\big(2V^2(t, x)I_2 - 2A_1^2(t, x)I_2 - \big(V^2(t, x)I_2 + A_1^2(t, x)I_2 -2A_1(t, x)V(t, x)\sigma_1\big)\nn\\
	&& - \big(V^2(t, x)I_2 + A_1^2(t, x)I_2 + 2A_1(t, x)V(t, x)\sigma_1\big)\big)\sigma_3\nn\\
	&=& -\big(-4A_1^2(t, x)I_2\big)\sigma_3 = 4A_1^2(t, x)\sigma_3.
	\eea
	In the derivation, we use the relations
	\be
	I_2\sigma_j = \sigma_jI_2, \quad j = 1, 3; \quad \sigma_1\sigma_3 = -\sigma_3\sigma_1.
	\ee
	Plugging \eqref{Ts1d1} and \eqref{Ts1d2} into \eqref{eq:TWd1}, we can obtain \eqref{eq:commu_1d} immediately.
	
	Similar derivation could be applied to the four-component Dirac equation \eqref{eq:dirac4} in 1D, and the details are omitted here for simplicity.\\

{\bf Appendix B}. Derivation of the double commutator in Lemma \ref{lemma:db_com_2d} for the Dirac equation \eqref{eq:dirac2} in 2D.
\setcounter{equation}{0}
\renewcommand{\theequation}{B.\arabic{equation}}

From \eqref{eq:split_TW}, in 2D, we have
\be
T = -\sigma_1\partial_1 - \sigma_2\partial_2 - i\sigma_3, \quad W(t) = -i\left(V(t, \bx)I_2 - A_1(t, \bx)\sigma_1 - A_2(t, \bx)\sigma_2\right).
\ee
Through the linearity of the double commutator in $T$,
\begin{equation}\label{eq:TWd2}
[W(t), [T, W(t)]] = -{\left[W(t), \left[\sigma_1\partial_1, W(t)\right]\right]} - [W(t), [\sigma_2\partial_2, W(t)]] - 
i{\left[W(t), \left[\sigma_3, W(t)\right]\right]}.
\end{equation}
From the definition of the Pauli matrices \eqref{eq:Paulim}, we have
\be\begin{split}
	\label{Pauliper}
	&\sigma_j^2=I_2, \quad \sigma_j\sigma_l=-\sigma_l\sigma_j, \qquad 1\le j\ne l\le 3,\\
	&\sigma_1\sigma_2=i\sigma_3, \quad \sigma_2\sigma_3=i\sigma_1,
	\quad \sigma_3\sigma_1=i\sigma_2.
\end{split}
\ee
Noticing \eqref{Pauliper}, we get
\begin{eqnarray}
	\label{2Dcommutator_1}
	&&[W(t), [\sigma_1\partial_1, W(t)]]\nn\\
	&&=  -\Big(2\big(V(t, \bx)I_2 - \sum_{j=1}^2A_j(t, \bx)\sigma_j\big)(\sigma_1\partial_1)
	\big(V(t, \bx)I_2 - \sum_{j=1}^2A_j(t, \bx)\sigma_j\big) \nn\\
	&&\ \ \ - \big(V(t, \bx)I_2 - \sum_{j=1}^2A_j(t, \bx)\sigma_j\big)^2(\sigma_1\partial_1) -
	(\sigma_1\partial_1)\big(V(t, \bx)I_2 - \sum_{j=1}^2A_j(t, \bx)\sigma_j\big)^2\Big)\nn\\
	&&=   -2\sigma_1A_2(t, \bx)\sigma_2\big(\partial_1V(t, \bx)I_2 - \sum_{j=1}^2A_j(t, \bx)\sigma_j\big)\nn\\
	&&\ \ \  - 2\sigma_1\big(V(t, \bx)I_2 - A_1(t, \bx)\sigma_1 + A_2(t, \bx)\sigma_2\big)
	\big(V(t, \bx)I_2 -\sum_{j=1}^2A_j(t, \bx)\sigma_j\big)\partial_1\nn\\
	&&\ \ \  + \sigma_1\big(V(t, \bx)I_2 - A_1(t, \bx)\sigma_1 +
	A_2(t, \bx)\sigma_2\big)^2\partial_1 + \sigma_1\big(V(t, \bx)I_2 - \sum_{j=1}^2A_j(t, \bx)\sigma_j\big)^2\partial_1\nn\\
	& &\ \ \ - 2\sigma_1A_2(t, \bx)\sigma_2
	\big(\partial_1V(t, \bx)I_2 - \sum_{j=1}^2A_j(t, \bx)\sigma_j\big)\nn\\
	&&=  -4A_2(t, \bx)\big(\partial_1V(t, \bx)\sigma_1\sigma_2 + \partial_1A_1(t, \bx)\sigma_2 -
	\partial_1A_2(t, \bx)\sigma_1\big) + 4A_2^2(t, \bx)\sigma_1\partial_1\nn\\
	&&\ \ \  - 4A_1(t, \bx)A_2(t, \bx)\sigma_2\partial_1\nn\\
	&&= 4\big(A_2^2(t, \bx)\sigma_1 - A_1(t, \bx)A_2(t, \bx)\sigma_2\big)\partial_1
	+ 4A_2(t, \bx)\big(\partial_1A_2(t, \bx)\sigma_1 - \partial_1A_1(t, \bx)\sigma_2\big)\\
	&&\ \ \ - 4iA_2(t, \bx)\partial_1V(t, \bx)\sigma_3\nn,
\end{eqnarray}
\begin{eqnarray}
	\label{2Dcommutator_3}
	[W(t), [\sigma_3, W(t)]] &= & -\Big(2\big(V(t, \bx)I_2 - \sum_{j=1}^2A_j(t, \bx)\sigma_j\big)
	\sigma_3\big(V(t, \bx)I_2 - \sum_{j=1}^2A_j(t, \bx)\sigma_j\big)\nn\\
	& &- \big(V(t, \bx)I_2 - \sum_{j=1}^2A_j(t, \bx)\sigma_j\big)^2\sigma_3
	- \sigma_3\big(V(t, \bx)I_2 - \sum_{j=1}^2A_j(t, \bx)\sigma_j\big)^2\Big)\nn\\
	&= & 2\sigma_3\big(V(t, \bx)I_2 + \sum_{j=1}^2A_j(t, \bx)\sigma_j\big)
	\sum_{j=1}^2A_j(t, \bx)\sigma_j \nn\\
	& &- 2\sigma_3\sum_{j=1}^2A_j(t, \bx)\sigma_j
	\big(V(t, \bx)I_2 - \sum_{j=1}^2A_j(t, \bx)\sigma_j\big)\nn\\
	&= & 4\big(A_1^2(t, \bx) + A_2^2(t, \bx)\big)\sigma_3,
\end{eqnarray}
and 
\begin{eqnarray}
	\label{2Dcommutator_2}
	[W(t), [\sigma_2\partial_2, W(t)]] &= &- 4\big(A_1(t, \bx)A_2(t, \bx)\sigma_1
	- A_1^2(t, \bx)\sigma_2\big)\partial_2 -4A_1(t, \bx)\big(\partial_2A_2(t, \bx)\sigma_1 - \partial_2A_1(t, \bx)\sigma_2\big)\nn\\
	&&+4iA_1(t, \bx)\partial_2V(t, \bx)\sigma_3 .
\end{eqnarray}
The derivation of \eqref{2Dcommutator_2} is similar to \eqref{2Dcommutator_1}, so the details are omitted for brevity. 
Plugging \eqref{2Dcommutator_1}, \eqref{2Dcommutator_3} and \eqref{2Dcommutator_2} into \eqref{eq:TWd2}, after some computation, we can get \eqref{eq:commu_2d}.

Similar derivation could be applied to the four-component Dirac equation \eqref{eq:dirac4} in 2D, and the details are omitted here for simplicity.\\

{\bf Appendix C}. Derivation of the double commutator in Lemma \ref{DTW3d} for the Dirac equation \eqref{eq:dirac4} in 3D.
\setcounter{equation}{0}
\renewcommand{\theequation}{C.\arabic{equation}}

The two operators $T$ and $W$ are defined as:
\be\label{TW3d1}
    T = -\sum_{j = 1}^3\alpha_j\partial_j -i\beta,\quad
    W(t) = -i\Bigl(V(t, \bx)I_4 - \sum_{j = 1}^3A_j(t, \bx)\alpha_j\Bigr).
\ee
By using the linearity of the double commutator in $T$, it is easy to obtain
   \be
	\begin{aligned}
	\label{split2}
	[W(t), [T, W(t)]] =& -[W(t), [\alpha_1\partial_1, W(t)]] - [W(t), [\alpha_2\partial_2, W(t)]]\\
	&-[W(t), [\alpha_3\partial_3, W(t)]] - i [W(t), [\beta, W(t)]].
	\end{aligned}
	\ee
	From \eqref{eq:alpha} and \eqref{gammam}, we have
	\be\begin{split}
		\label{matricper}
		&\beta^2=I_4, \quad \alpha_j^2=I_4, \quad \alpha_j\alpha_l=-\alpha_l\alpha_j,\\
		&\beta\alpha_j = -\alpha_j\beta,
		\quad \gamma\alpha_j=\alpha_j\gamma,  \qquad 1\le j\ne l\le 3,\\
		&\alpha_1\alpha_2=i\gamma\alpha_3, \quad \alpha_2\alpha_3=i\gamma\alpha_1, \quad \alpha_3\alpha_1=i\gamma\alpha_2.
	\end{split}
	\ee
	Noticing \eqref{TW3d1}, and \eqref{matricper}, we get
	\begin{eqnarray}
	\label{3Dcommutator_4}
	[W(t), [\beta, W(t)]] &= & -\bigg(2\Big(V(t)I_4 - \sum_{j = 1}^3A_j(t)\alpha_j\Big)\beta
	\Big(V(t)I_4 - \sum_{j = 1}^3A_j(t)\alpha_j\Big)\nn\\
	&& - \Big(V(t)I_4 - \sum_{j = 1}^3A_j(t)\alpha_j\Big)^2\beta
	- \beta\Big(V(t)I_4 - \sum_{j = 1}^3A_j(t)\alpha_j\Big)^2\bigg)\nn\\
	&= & -2\beta\Big(V(t)I_4 + \sum_{j = 1}^3A_j(t)\alpha_j\Big)
	\Big(V(t)I_4 - \sum_{j = 1}^3A_j(t)\alpha_j\Big)\nn\\
	&& + \beta\Big(V(t)I_4 + \sum_{j = 1}^3A_j(t)\alpha_j\Big)^2 +
	\beta\Big(V(t)I_4 - \sum_{j = 1}^3A_j(t)\alpha_j\Big)^2\nn\\
	&= & 4\big(A_1^2(t) + A_2^2(t) + A_3^2(t)\big)\beta.
	\end{eqnarray}
	\bea\label{3Dcommutator_1}
	&&[W(t), [\alpha_1\partial_1, W(t)]]\nn\\
	&&=  -\bigg(2\Big(V(t)I_4
	- \sum_{j = 1}^3A_j(t)\alpha_j\Big)(\alpha_1\partial_1)\Big(V(t)I_4
	- \sum_{j = 1}^3A_j(t)\alpha_j\Big)\nn\\
	&&\ \ \ - \Big(V(t)I_4 - \sum_{j = 1}^3A_j(t)\alpha_j\Big)^2(\alpha_1\partial_1)
	- (\alpha_1\partial_1)\Big(V(t)I_4 - \sum_{j = 1}^3A_j(t)\alpha_j\Big)^2\bigg)\nn\\
	&&=  -4\alpha_1\big(A_2(t)\alpha_2 + A_3(t)\alpha_3\big)\big(\partial_1V(t)I_4 -
	\partial_1A_1(t)\alpha_1 - \partial_1A_2(t)\alpha_2 - \partial_1A_3(t)\alpha_3\big)\nn\\
	&&\ \  \ + \alpha_1\bigg(\Big(V(t)I_4 - A_1(t)\alpha_1 + A_2(t)\alpha_2
	+ A_3(t)\alpha_3\Big)^2
	+ \Big(V(t)I_4 - \sum_{j = 1}^3A_j(t)\alpha_j\Big)^2\nn\\
	&&\ \ \  - 2\Big(V(t)I_4 - A_1(t)\alpha_1 +
	A_2(t)\alpha_2 + A_3(t)\alpha_3\Big)\Big(V(t)I_4 - \sum_{j = 1}^3A_j(t)\alpha_j\Big)\bigg)\partial_1,\nn\\
	&&=  4\big(A_2(t)\alpha_2 + A_3(t)\alpha_3\big)\alpha_1\big(\partial_1V(t)I_4 -
	\partial_1A_1(t)\alpha_1 - \partial_1A_2(t)\alpha_2 - \partial_1A_3(t)\alpha_3\big)\nn\\
	&&\ \ \ + 4\Big(\big(A_2^2(t) + A_3^2(t)\big)\alpha_1 - A_1(t)A_2(t)\alpha_2 -
	A_1(t)A_3(t)\alpha_3\Big)\partial_1\nn\\
	&&= 4\Big(\big(A_2(t)\partial_1A_2(t) + A_3(t)\partial_1A_3(t)\big)\alpha_1
	- A_2(t)\partial_1A_1(t)\alpha_2 -  A_3(t)\partial_1A_1(t)\alpha_3\nn\\
	&&\ \ \ + \big(iA_2(t)\partial_1A_3(t) - iA_3(t)\partial_1A_2(t)\big)\gamma
	+ iA_3(t)\partial_1V(t)\gamma\alpha_2 - iA_2(t)\partial_1V(t)\gamma\alpha_3\Big)\nn\\
	&&\ \ \ + 4\Big(\big(A_2^2(t) + A_3^2(t)\big)\alpha_1 - A_1(t)A_2(t)\alpha_2 -
	A_1(t)A_3(t)\alpha_3\Big)\partial_1.
\end{eqnarray}
\begin{eqnarray}
\label{3Dcommutator_2}
&&[W(t), [\alpha_2\partial_2, W(t)]]\nn\\
&&= 4\Big(-A_1(t)\partial_2A_2(t)\alpha_1
+ \big(A_1(t)\partial_2A_1(t) + A_3(t)\partial_2A_3(t)\big)\alpha_2 - A_3(t)\partial_2A_2(t)\alpha_3\nn\\
&&\ \ \  + \big(iA_3(t)\partial_2A_1(t) - iA_1(t)\partial_2A_3(t)\big)\gamma
- iA_3(t)\partial_2V(t)\gamma\alpha_1 + iA_1(t)\partial_2V(t)\gamma\alpha_3\Big)\nn\\
&&\ \ \  + 4\Big(\big(A_1^2(t) + A_3^2(t)\big)\alpha_2 - A_2(t)A_1(t)\alpha_1 -
A_2(t)A_3(t)\alpha_3\Big)\partial_2.
\end{eqnarray}
\begin{eqnarray}
\label{3Dcommutator_3}
&&[W(t), [\alpha_3\partial_3, W(t)]]\nn\\
&&= 4\Big(- A_1(t)\partial_3A_3(t)\alpha_1
- A_2(t)\partial_3A_3(t)\alpha_2 + \big(A_1(t)\partial_3A_1(t) + A_2(t)\partial_3A_2(t)\big)\alpha_3\nn\\
& &\ \ \ + \big(iA_1(t)\partial_3A_2(t) - iA_2(t)\partial_3A_1(t)\big)\gamma +
iA_2(t)\partial_3V(t)\gamma\alpha_1 - iA_1(t)\partial_3V(t)\gamma\alpha_2\Big)\nn\\
&&\ \ \ + 4\Big(\big(A_1^2(t) + A_2^2(t)\big)\alpha_3 - A_3(t)A_1(t)\alpha_1 -
A_3(t)A_2(t)\alpha_2\Big)\partial_3.
\end{eqnarray}
In the above, we use $V(t) := V(t, \bx)$ and $A_j(t) := A_j(t, \bx)$, $j = 1, 2, 3$, for brevity.

Plugging (\ref{3Dcommutator_1}), (\ref{3Dcommutator_2}), (\ref{3Dcommutator_3})
and (\ref{3Dcommutator_4}) into (\ref{split2}), after some computation,
we could obtain \eqref{commun_3d}.

\end{document}